\documentclass{article}%
\usepackage{amssymb}
\usepackage{amsmath}
\usepackage{graphicx}
\usepackage{amsfonts}%
\providecommand{\U}[1]{\protect\rule{.1in}{.1in}}
\begin{document}

\title{Determination through Universals:\\An Application of Category Theory in the Life Sciences}
\author{David Ellerman\\Philosophy Department\\University of California at Riverside}
\date{}
\maketitle

\begin{abstract}
Category theory has foundational importance because it provides conceptual
lenses to characterize what is important and universal in mathematics---with
adjunctions being the primary lense. If adjunctions are so important in
mathematics, then perhaps they will isolate concepts of some importance in the
empirical sciences. But the applications of adjunctions have been hampered by
an overly restrictive formulation that avoids heteromorphisms or hets. By
reformulating an adjunction using hets, it is split into two parts, a left and
a right semiadjunction (or half-adjunuction). Semiadjunctions (essentially a
formulation of a universal mapping property using hets) turn out to be the
appropriate concept for applications in the life sciences. The semiadjunctions
characterize three principal schemes with applications:

1. determination through a receiving universal (e.g., natural selection,
perception, language acquisition, recursion, and language understanding);

2. determination through a sending universal (e.g., intentional action, DNA
mechanism to build amino acids, hierarchy of regulatory genes to build organs,
stem cells, and language production), and

3. two-way determination with one universal as both receiving and sending
(e.g., perception/action and language understanding/production).

\end{abstract}
\tableofcontents

\section{Introduction}

Category theory has foundational importance because it provides conceptual
lenses to characterize what is important and universal in mathematics---with
an adjunction (or pair of adjoint functors) being the primary lense. The
mathematical importance of adjunctions is now well recognized. As Steven
Awodey put it in his recent text:

\begin{quote}
{\footnotesize The notion of adjoint functor applies everything that we have
learned up to now to unify and subsume all the different universal mapping
properties that we have encountered, from free groups to limits to
exponentials. But more importantly, it also captures an important mathematical
phenomenon that is invisible without the lens of category theory. Indeed, I
will make the admittedly provocative claim that adjointness is a concept of
fundamental logical and mathematical importance that is not captured elsewhere
in mathematics. }$\cite[p. 179]{aw:cat}$
\end{quote}

\noindent Other category theorists have given similar testimonials.

\begin{quote}
{\footnotesize To some, including this writer, adjunction is the most
important concept in category theory. \cite[p. 6]{wood:ord}}

{\footnotesize The isolation and explication of the notion of
\textit{adjointness} is perhaps the most profound contribution that category
theory has made to the history of general mathematical ideas.\cite[p.
438]{Gold:topoi}}

{\footnotesize Nowadays, every user of category theory agrees that
[adjunction] is the concept which justifies the fundamental position of the
subject in mathematics. \cite[p. 367]{tay:pfm}}
\end{quote}

If a concept, like that of a pair of adjoint functors, is of such importance
in mathematics, then one might expect it to have applications, perhaps of some
importance, in the empirical sciences. Yet this does not seem to be the case,
particularly in the life sciences. Perhaps the problem has been finding the
right level of generality or specificity where non-trivial applications can be
found, i.e., finding out "where theory lives."

This paper argues that the application of adjoints has been hampered by an
overly specific formulation of the adjunctive properties that only uses
homomorphisms or homs\footnote{"Hom" is pronounced to rhyme with "Tom" or
"bomb."} (object-to-object morphisms \textit{within} a category). A
reformulation of adjunctions using heteromorphisms or hets (object-to-object
morphisms between objects of \textit{different} categories) allows an
adjunction to be split into two "semiadjunctions." The argument is that a
semiadjunction (essentially a reformulation of a universal mapping property
using hets) turns out to be the right concept for applications. Moreover, by
"splitting the atom" of an adjunction into two semiadjunctions, the
semiadjunctions can be recombined in a different way to define the cognate
notion of a "brain functor"--which, as the name indicates, has applications in
cognitive science.

This application of category theory in the life sciences is not exact in the
sense that, say, calculus is applied in physics. The category theory is exact
at the mathematical level, but it is the general determinative schema of the
semiadjunctions that is applied.

There is already a considerable but widely varying literature on the
application of category theory to the life sciences--such as the work of
Robert Rosen \cite{rosen:life} and his followers\footnote{See
\cite{zafiris:rosen} and that paper's references.} as well as Andr\'{e}e
Ehresmann and Jean-Paul Vanbremeersch \cite{ehres:mes} and their
commentators.\footnote{See \cite{kainen:ehres} for Kainen's comments on the
Ehresmann-Vanbremeersch approach, Kainen's own approach, and a broad
bibliography of relevant papers.} But it is still early days, and many
approaches need to be tried to find out "where theory lives."

The approach taken here is based on a specific use of the characteristic
concepts of category theory, namely universal mapping properties, to define a
general schema of determination through universals. The closest approach in
the literature (but without the hets and semiadjunctions) is that of
Fran\c{c}ois Magnan and Gonzalo Reyes \cite{mag-reyes:cog} which emphasizes
that "Category theory provides means to circumscribe and study what is
universal in mathematics and other scientific disciplines." \cite[p.
57]{mag-reyes:cog}. Their intended field of application is cognitive science.

\begin{quotation}
We may even suggest that universals of the mind may be expressed by means of
universal properties in the theory of categories and much of the work done up
to now in this area seems to bear out this suggestion....

By discussing the process of counting in some detail, we give evidence that
this universal ability of the human mind may be conveniently conceptualized in
terms of this theory of universals which is category theory. \cite[p.
59]{mag-reyes:cog}\footnote{See the example of recursion developed below.}
\end{quotation}

The approach of determination through universals, like any approach, should be
judged on how well it isolates and describes the essential and important
features of biological and cognitive systems.

\section{Adjunctions, hets, and semiadjunctions}

In the body of this paper, I will try to keep the mathematics at a minimal
conceptual level--which the mathematical formulations restricted to the
Appendix. Category theory lends itself to visualization in diagrams so that
non-mathematical style of presentation is emphasized by an abundant use of diagrams.

A \textit{category} is intuitively a set of objects of the same type.
\textit{Morphisms} between objects should be thought of as a type of
determining relation or cause-effect relation between the objects. When a
morphism is between objects of the same category, it is called a
\textit{homomorphism} or \textit{hom}, and when between objects of different
categories it is a \textit{heteromorphism} or \textit{het}.\footnote{See
\cite{ell:taf} or \cite{ell:ae} for the introduction of heteromorphisms in
category theory.} One of the problems in the conventional treatment of
category theory\footnote{The standard presentation is Mac Lane's text
\cite{mac:cwm} but Lawvere and Schanuel's book \cite{law:conceptual} is a more
conceptual introduction. Magnan and Reyes \cite{mag-reyes:cog} give an
excellent informal treatment of the main universal constructions.} is that it
tries to ignore heteromorphisms even though hets are a natural part of working
mathematics. This leads to certain definitions being rather contrived (to
avoid mentioning hets), the usual treatment of adjunctions being the case in point.

Adjunctions will be introduced informally and in the natural manner using
hets. The general setting is how the objects in one category (e.g., the
"environment" in a life sciences context), the "sending" category, will
"affect" or "determine" objects in another category (e.g., "organisms"), the
"receiving" category. We start with an object $X$ in the sending category, an
object $X$ in the receiving category, and a specific het determination
$d:X\dashrightarrow A$ from $X$ to $A$.\footnote{Hets are represented by
dashed arrows $\dashrightarrow$ and homs by solid arrows $\rightarrow$.}%

\begin{center}
\includegraphics[
height=2.0897in,
width=4.3088in
]%
{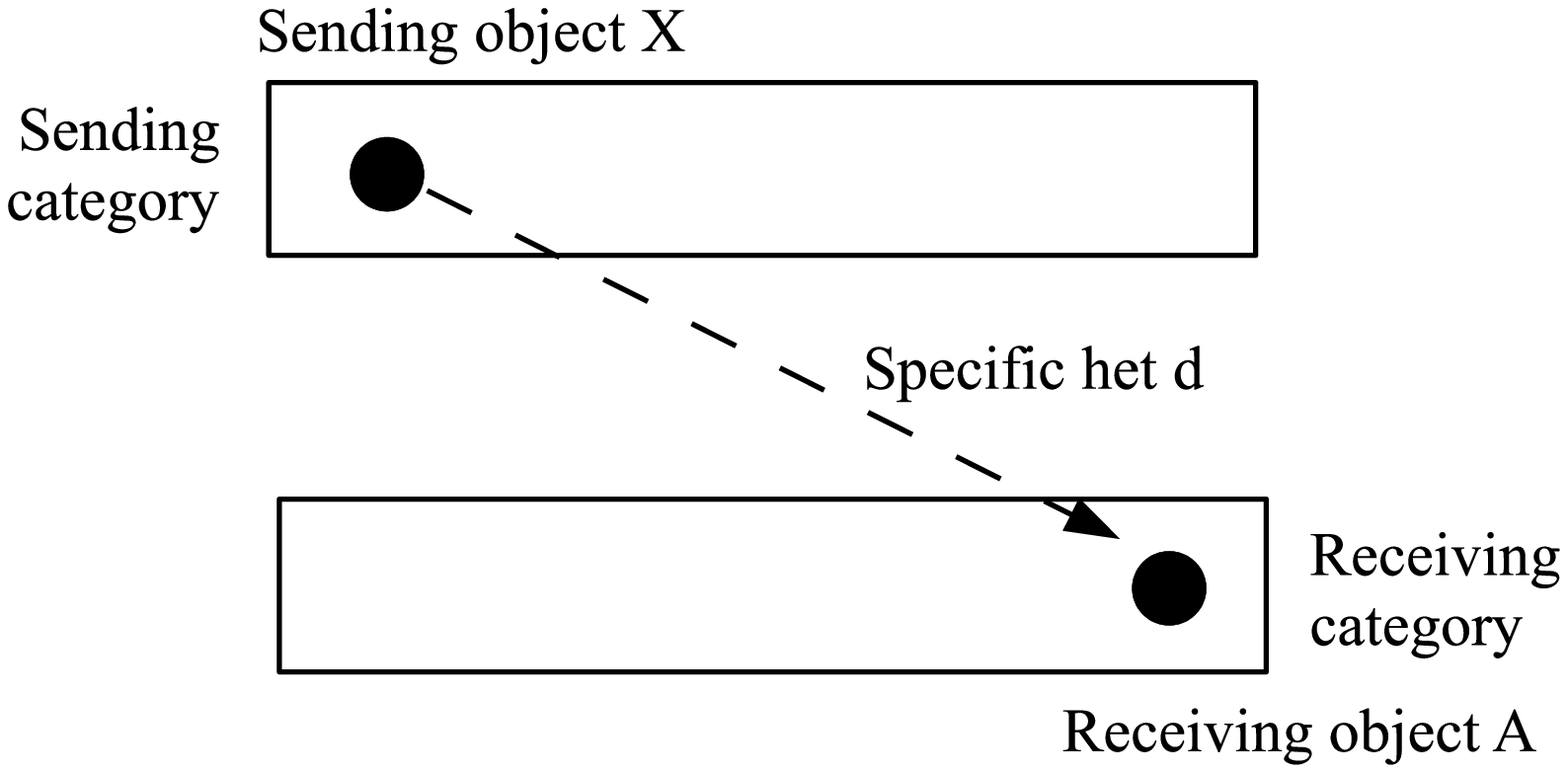}%
\\
Figure 1: Het $d:X\dashrightarrow A$%
\end{center}

One of the most important concepts isolated by category theory (more basic
than adjunctions) is the concept of a \textit{universal mapping property} or
UMP\footnote{A UMP illustrates the philosophical notion of a concrete
universal \cite{ell:cu}.} which models an important type of determination,
\textit{determination through universals}. With each sending object $X$ in the
sending category, there is an associated \textit{receiving universal} $F(X)$
in the receiving category, and there is a \textit{universal receiving het}
$h_{X}:X\dashrightarrow F(X)$. The universal mapping property is: for every
het $d:X\dashrightarrow A$, there is a unique hom $f(d):F(X)\longrightarrow A$
in the receiving category such that:

\begin{center}
$X\overset{h_{X}}{\dashrightarrow}F(X)\overset{f(d)}{\longrightarrow
}A=X\overset{d}{\dashrightarrow}A$
\end{center}

\noindent i.e., such that the determination through the receiving universal
het $h_{X}:X\dashrightarrow F(X)$ followed by the hom $f(d):F(X)\rightarrow A$
is the same as the original het $d:X\dashrightarrow A$.%

\begin{center}
\includegraphics[
height=2.0432in,
width=4.43in
]%
{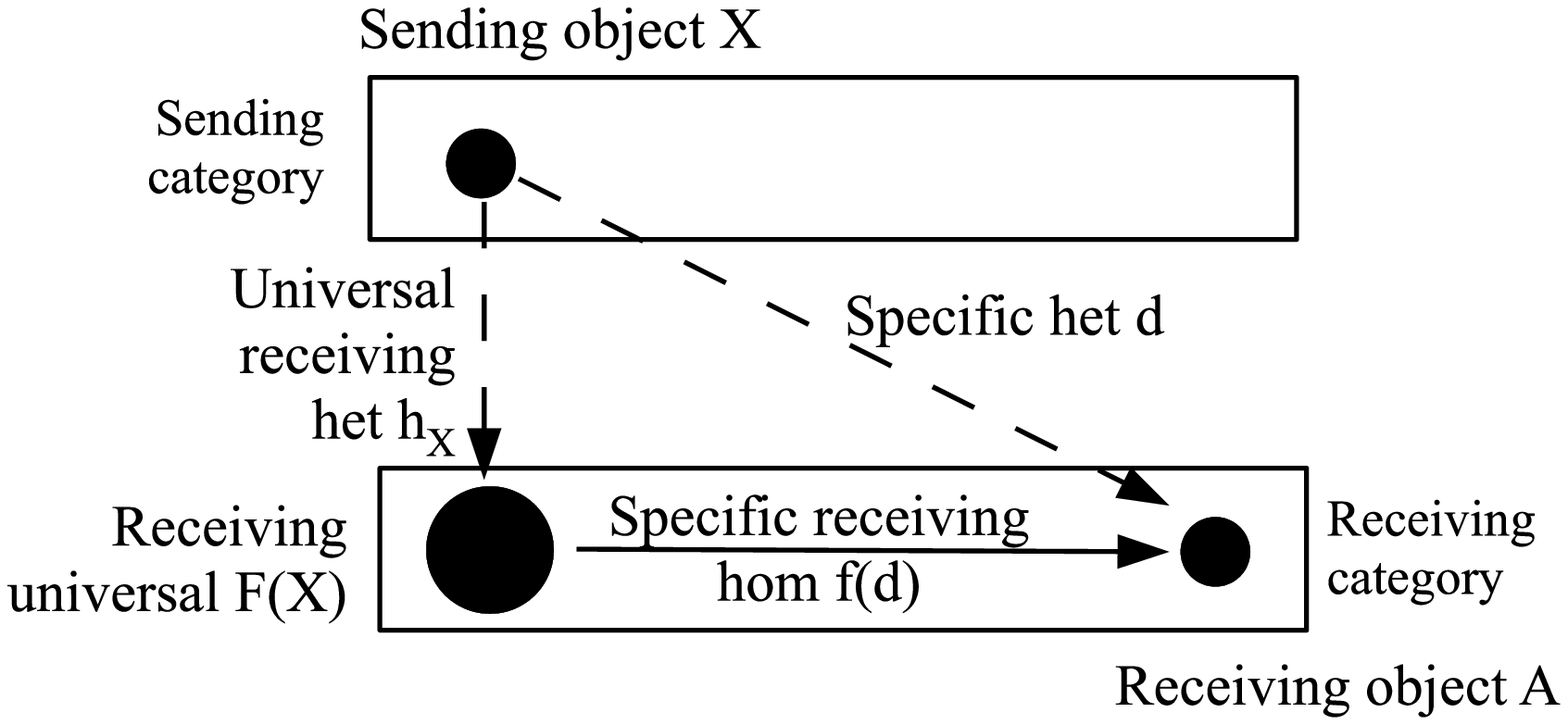}%
\\
{}Figure 2: Scheme for determination by a receiving universal $F(X)$.
\end{center}

Moreover given any hom $f:F(X)\rightarrow A$ in the receiving category,
preceding it by the universal receiving het $h_{X}$ would give a het
$fh_{X}:X\dashrightarrow A$ which would play the role of $d:X\dashrightarrow
A$ in the above equation. Hence there is a one-to-one correspondence (or
isomorphism) between the hets $d:X\dashrightarrow A$ and the homs
$f:F(X)\rightarrow A$. This isomorphism is also canonical or natural (in the
category-theoretic sense), and can be represented as:

\begin{center}
$\operatorname*{Hom}_{receiving}(F(X),A)\cong\operatorname*{Het}(X,A)$.
\end{center}

\noindent\ \noindent Even before completing our definition of an adjunction,
we can define the above situation, given by the association of the receiving
universal $F(X)$ with each object $X$ in the sending category (and similarly
for morphisms so it is a "functor") along with the canonical isomorphism
$\operatorname*{Hom}(F(X),A)\cong\operatorname*{Het}(X,A)$, as a \textit{left
semiadjunction. }Even if a full adjunction is not present, we will still refer
to the functor $F()$ as a \textit{left adjoint} whose values $F(X)$ are the
receiving universals.

Then we return to the specific het $d$ and take the dual case which will
define the dual notion of a "right semiadjunction." With each receiving object
$A$ in the receiving category, there is an associated \textit{sending
universal} $G(A)$ in the sending category, and there is a \textit{sending
universal het} $e_{A}:G(A)\dashrightarrow A$. The universal mapping property
is: for every het $d:X\dashrightarrow A$, there is a unique hom
$g(d):X\longrightarrow G(A)$ in the sending category such that:

\begin{center}
$X\overset{g(d)}{\longrightarrow}G(A)\overset{e_{A}}{\dashrightarrow
}A=X\overset{d}{\dashrightarrow}A$
\end{center}

\noindent i.e., such that the determination through the universal sending het
$e_{A}:G(A)\dashrightarrow A$ preceded by the hom $g(d):X\rightarrow G(A)$ is
the same as the original het $d:X\dashrightarrow A$.%

\begin{center}
\includegraphics[
height=2.054in,
width=4.4109in
]%
{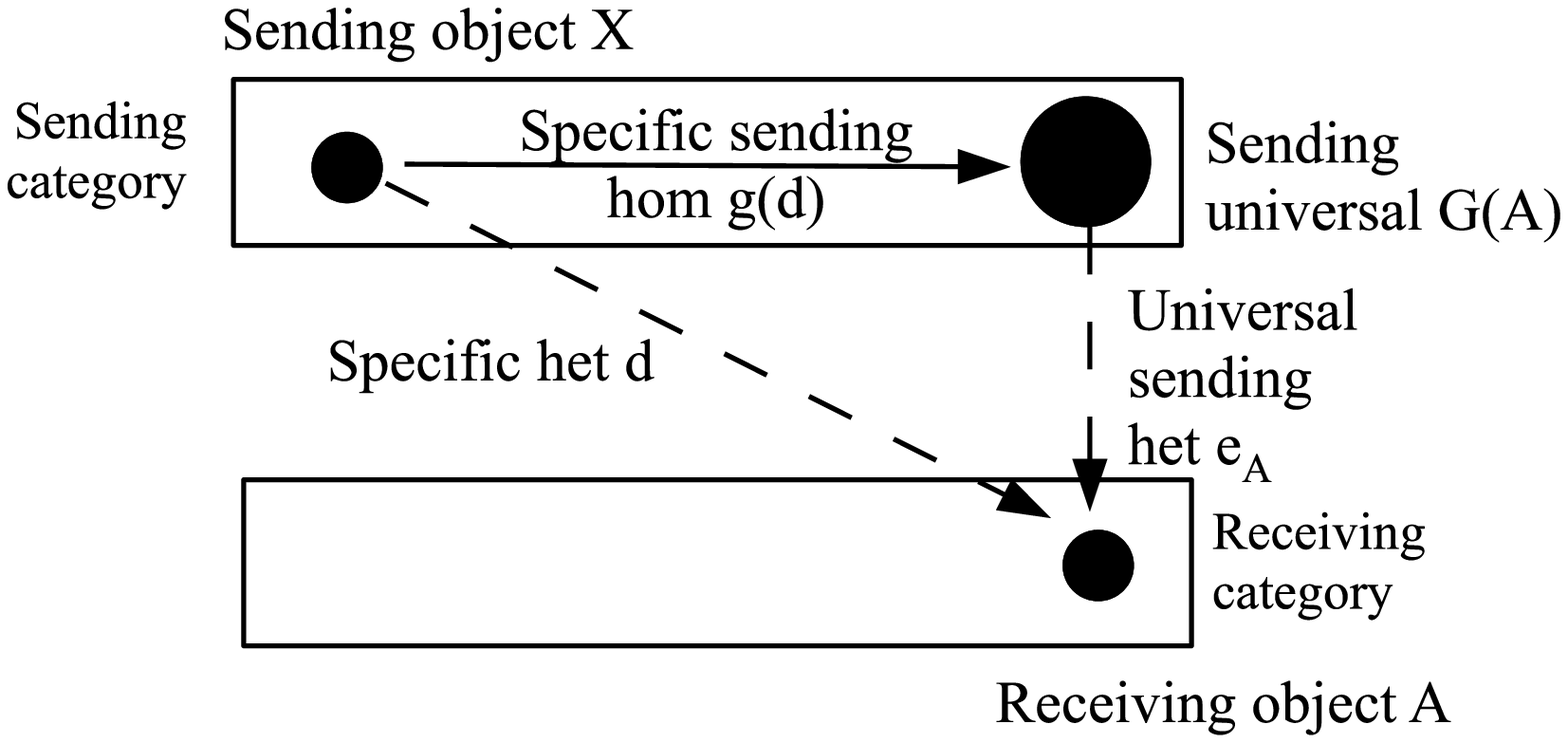}%
\\
Figure 3: Scheme for determination by a sending universal $G(A)$.
\end{center}

Moreover given any hom $g:X\rightarrow G(A)$ in the sending category,
following it by the universal sending het $e_{A}$ would give a het
$e_{A}g:X\dashrightarrow A$ which would play the role of $d:X\dashrightarrow
A$ in the above equation. Hence there is a one-to-one correspondence (or
isomorphism) between the hets $d:X\dashrightarrow A$ and the homs
$g:X\rightarrow G(A)$. This isomorphism is also canonical or natural, and can
be represented as:

\begin{center}
$\operatorname*{Het}(X,A)\cong\operatorname*{Hom}_{sending}(X,G(A))$.
\end{center}

Dually, we can define the above situation, given by the association of the
sending universal $G(A)$ with each object $A$ in the receiving category along
with the canonical isomorphism $\operatorname*{Het}(X,A)\cong%
\operatorname*{Hom}_{sending}(X,G(A))$, as a \textit{right semiadjunction.
}Even if a full adjunction is not present, we will still refer to the functor
$G()$ as a \textit{right adjoint} whose values are the sending universal
$G(A)$.

Now we are prepared to define an \textit{adjunction} essentially as:

\begin{center}
adjunction = left semiadjunction + right semiadjunction

$\operatorname*{Hom}_{receiving}(F(X),A)\cong\operatorname*{Het}%
(X,A)\cong\operatorname*{Hom}_{sending}(X,G(A))$.
\end{center}

%

\begin{center}
\includegraphics[
height=2.3769in,
width=4.8144in
]%
{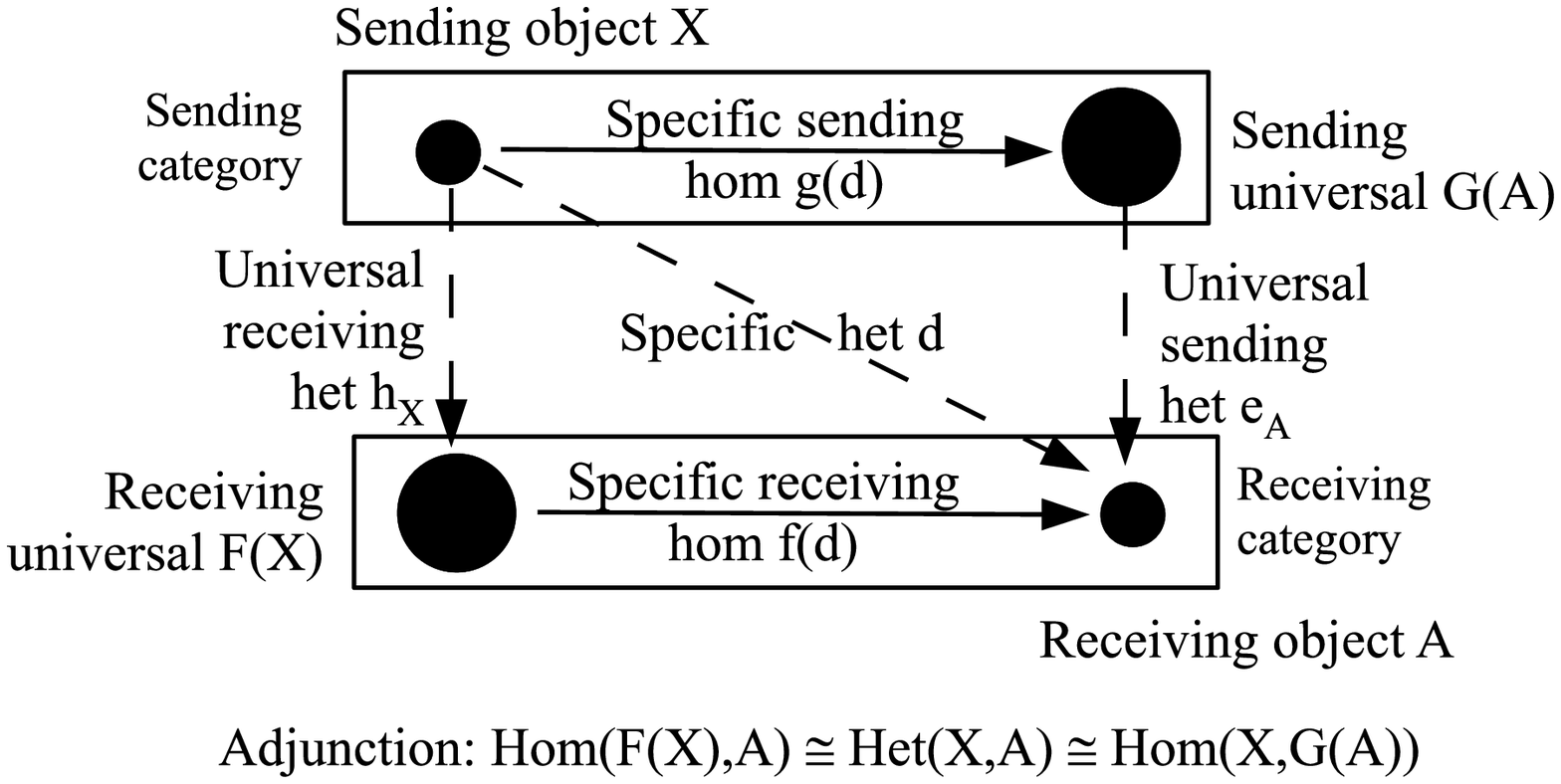}%
\\
Figure 4: The Adjunctive Square Diagram
\end{center}

Recall that the conventional treatment of an adjunction keeps the hets "in the
closet," so the defining natural isomorphism just leaves out the middle het
term $\operatorname*{Het}(X,A)$ to give the het-free notion of an adjunction:
$\operatorname*{Hom}(F(X),A)\cong\operatorname*{Hom}(X,G(A))$ where $F()$ is
the \textit{left adjoint} and $G()$ is the \textit{right adjoint}.

In almost all adjunctions, one of the adjoints is rather trivial and contrived
so the principal universal mapping property is expressed by the other adjoint,
but one needs the trivial adjoint in order to avoid the hets. This has
hampered applications since it is the main receiving or sending universal that
is important, not the trivial other adjoint needed to formulate the het-free
adjunction. Moreover since determinations between quite different type of
objects are important in applications, hets are an essential part of the
story, and thus the treatment of adjunction leaving out the hets also hampers
applications. But by bringing the hets "out of the closet," we can "fission"
the adjunction into two semiadjunctions, one of which is typically appropriate
in an application. And we can recombine the two semiadjunctions in a different
way to define a "brain functor"--which combines two non-trivial dual
semiadjunctions in one scheme of \textit{two-way} determination.

\section{Determinations through receiving universals}

\subsection{Selection versus instruction}

In the most generic example of determination through a receiving universal
(i.e., the determinative scheme of a left semiadjunction), the sending
category or domain is the "environment" while the receiving category or domain
is an "organism" (or population of organisms) and the factoring of a het
determination $d:X\dashrightarrow A$ through the receiving universal $F(X)$ is
a type of selective "recognition." Indeed, the contrast between the direct
determination by the het $d:X\dashrightarrow A$ and the factorization through
the receiving universal can be seen as the contrast between an instructionist
account and a selectionist account of some determination.

We begin with a toy example from Peter Medawar \cite[pp. 88-90]{med:fut}. A
jukebox contains a "universal" set of records each of which has the
instructions for some piece of music. A person in the environment selects and
pushes a button which causes the jukebox to supply the musical instructions in
the form of a record to the record player part of the mechanism.
Alternatively, the person in the environment could directly supply the musical
instructions in the form of a record to a record player. In each case, we end
up with a record player playing a record but by two different types of
mechanisms, one being selective and the other being instructive.%

\begin{center}
\includegraphics[
height=1.9767in,
width=5.0568in
]%
{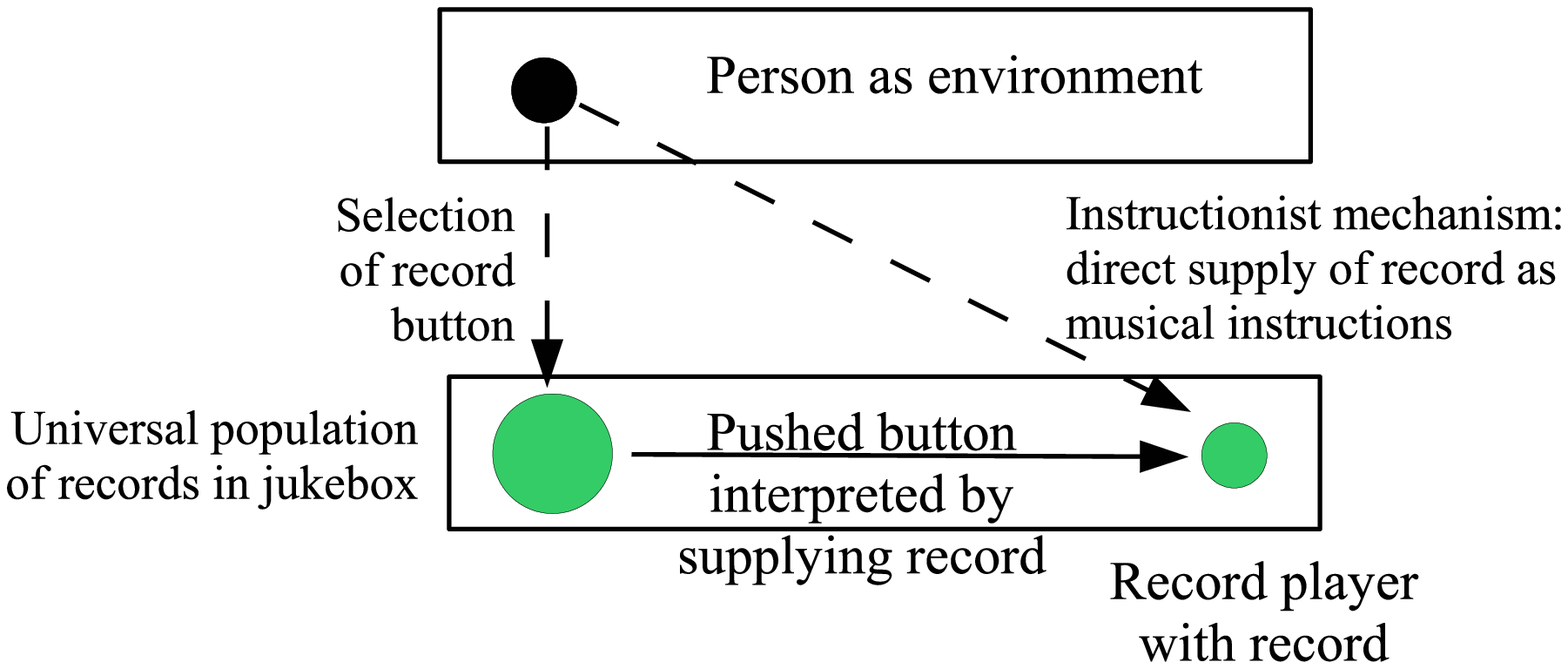}%
\\
Figure 5: Selection as determination through a receiving universal
\end{center}

\subsection{Evolution by natural selection}

After Gerald Edelman received the Nobel prize for his work on the selectionist
approach to the immune system, he switched to neurophysiology and is thus
well-placed to isolate the selectionist commonalities in these different
domains (which turn out to be modeled by left semiadjunctions).

\begin{quote}
{\footnotesize The long trail from antibodies to conscious brain events has
reinforced my conviction that evolution, immunology, embryology, and
neurobiology are all sciences of recognition whose mechanics follow
selectional principles. ...All selectional systems follow three principles.
There must be a generator of diversity, a polling process across the diverse
repertoires that ensue, and a means of differential amplification of the
selected variants.\cite[p. 7367]{edel:recog}; (also \cite[pp.41-2]{edel:sky})
}
\end{quote}

\noindent These three principles are functionally represented by the three
components in a determination through a receiving universal. The
\textquotedblleft generator of diversity\textquotedblright\ is the receiving
universal object, the \textquotedblleft polling process across the diverse
repertoires\textquotedblright\ is represented by the receiving universal
morphism that is the canonical external-internal interface between external
environment and the receiving universal object, and finally the
\textquotedblleft differential amplification\textquotedblright\ is represented
by the factor hom.%

\begin{center}
\includegraphics[
height=2.6351in,
width=5.0153in
]%
{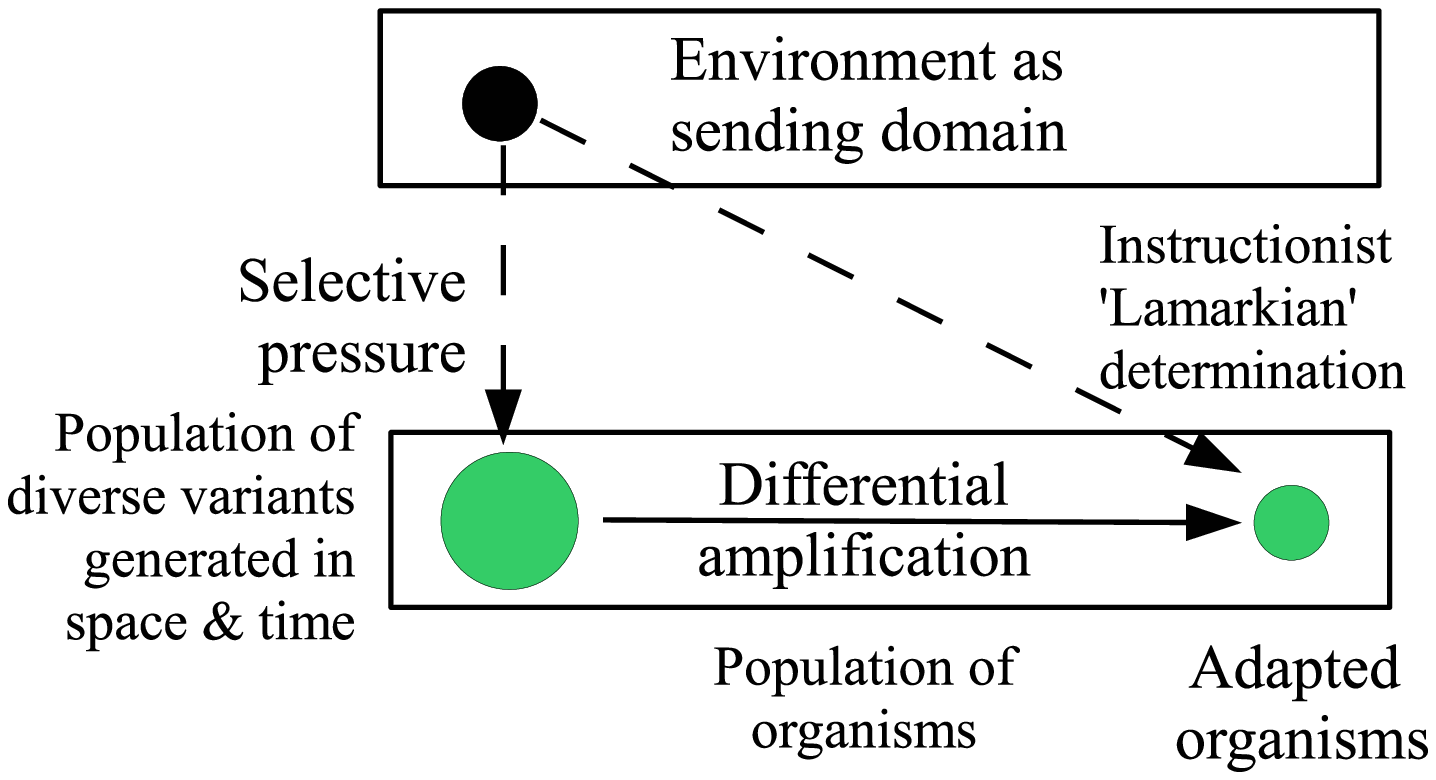}%
\\
Figure 6: Natural selection as determination through a receiving universal.
\end{center}

The direct determination given by the het would represent a Lamarckian
instructionist determination; the environment directly inducing the changes in
the organism to adapt it to the environment. The indirect determination
through the receiving universal using the "universal" variations in organisms
as the basis for the environment's Darwinian natural selection so that the
adapted organisms differentially reproduce.

\subsection{Immune system as a selectionist mechanism}

Evolution is not the only example in the life sciences of determinative
processes that were originally assumed to be instructionist but were later
found to operate by a selectionist mechanism. One of the most telling cases
was the immune system. Originally it was assumed that the antigen would
somehow instruct the immune system as to how an anti-body could be constructed
to neutralize the antigen. During the 1950s, a number of difficulties in the
instructionist account fostered the development of a selectionist approach.
While many researchers contributed to this approach, one of the earliest was
Niels Jerne \cite{jerne:nst} who has also been most attentive to analogies
with other fields.

In the selectionist theory, the immune system takes on the active role of
generating a huge well-nigh \textquotedblleft universal\textquotedblright%
\ variety of antibodies but in low concentrations. This initial generation of
candidate antibodies is not being directed or instructed by the past disease
history of the organism. An externally introduced antigen has the indirect
role of simply selecting which antibody fits it like a key in a lock. Every
antibody has the possibility of self-reproducing or cloning itself but it is
the ones whose key has fit into a lock that have this potentiality triggered.
Then that antibody is differentially amplified in the sense of being cloned
into many copies to lock up the other instances of the antigen. Thus the
selectionist account of the immune system has the main features of
determination through a receiving universal.%

\begin{center}
\includegraphics[
height=2.4583in,
width=4.792in
]%
{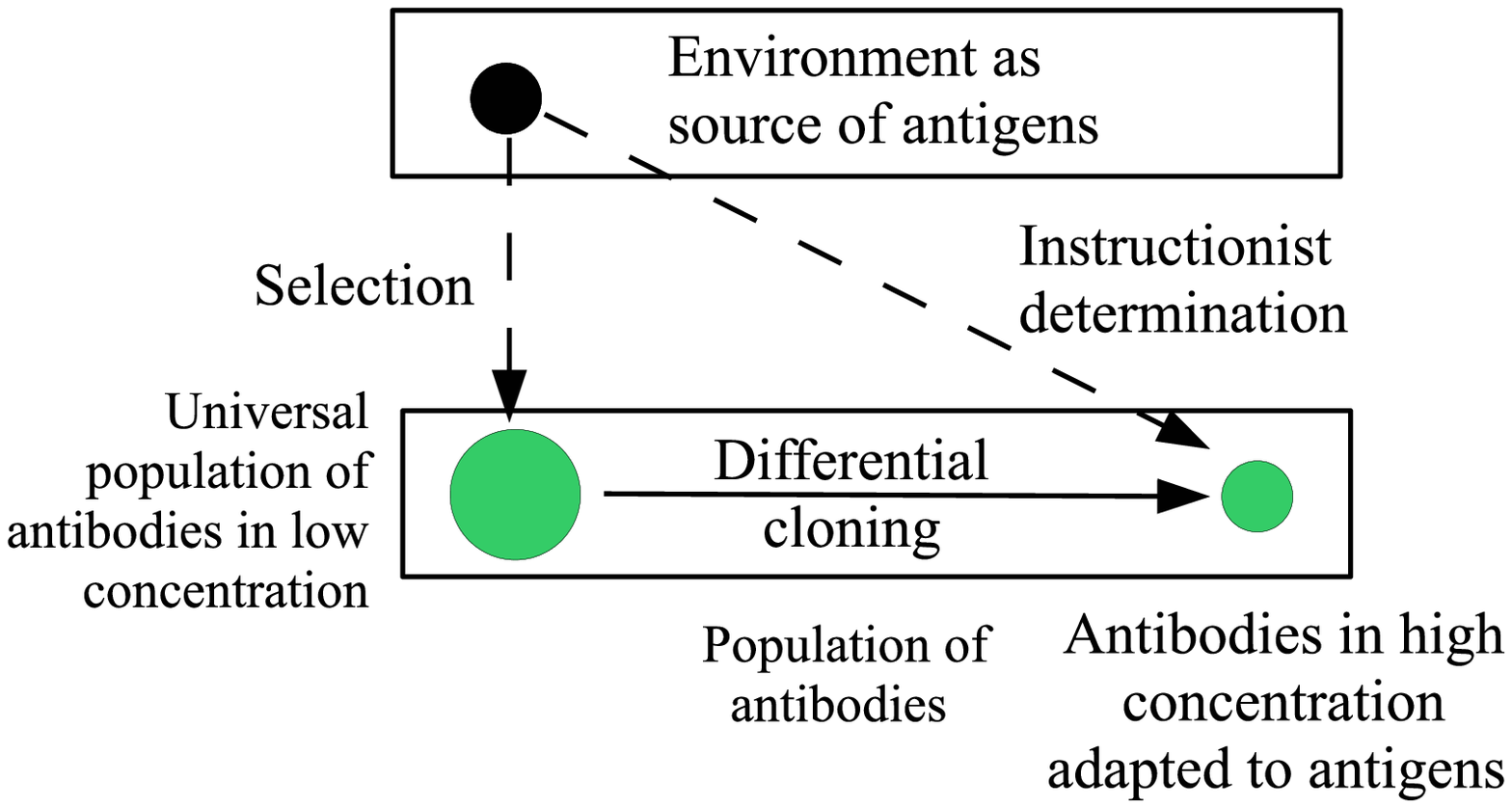}%
\\
Figure 7: Immune system as determination through a receiving universal.
\end{center}

\subsection{Edelman's neural Darwinism}

Edelman carried over the selectionist model to neurophysiology to develop his
theory of neuronal group selection or neural Darwinism.

\begin{quote}
{\footnotesize [T]he theoretical principle I shall elaborate here is that the
origin of categories in higher brain function is somatic selection among huge
numbers of variants of neural circuits contained in networks created
epigenetically in each individual during its development; this selection
results in differential amplification of populations of synapses in the
selected variants. In other words, I shall take the view that the brain is a
selective system more akin in its workings to evolution than to computation or
information processing.\cite[p. 25]{edel:nd}}
\end{quote}

\noindent There are several different phases in this selectionist theory. In
the developmental phase of the brain, a huge variety of loose connections are
made. Those that find some resonance with the individual's experience are
strengthen while those that are unused will atrophy. The slogan is that
\textquotedblleft the neurons that fire together, wire
together.\textquotedblright\ Later there is an experiential selection the
strengthens some connections and weakens others. Finally, \textquotedblleft
reentrant\textquotedblright\ signals within the brain deepen the process of
self-organization through strengthening some connections and weakening others.

In broad-brush terms, one might intuitively think of the universal model as a
large set of brain circuits representing a wide (\textquotedblleft
universal\textquotedblright) range of sensory images and vibrating at a low
level of amplitude beneath the level of consciousness (analogous to the
\textquotedblleft universal\textquotedblright\ repertoire of antibodies
present in the immune system in low concentrations). When a specific signal is
received from the environment, then it might resonate with a particular
circuit-image which would greatly increase the amplitude of those vibrations
and would thus constitute the perception. This sort of model has a type of
intentionality (i.e., seeing is always seeing-as) since the perception would
always be \textquotedblleft perception-as\textquotedblright\ depending on
which image was resonated.%

\begin{center}
\includegraphics[
height=2.6418in,
width=5.3001in
]%
{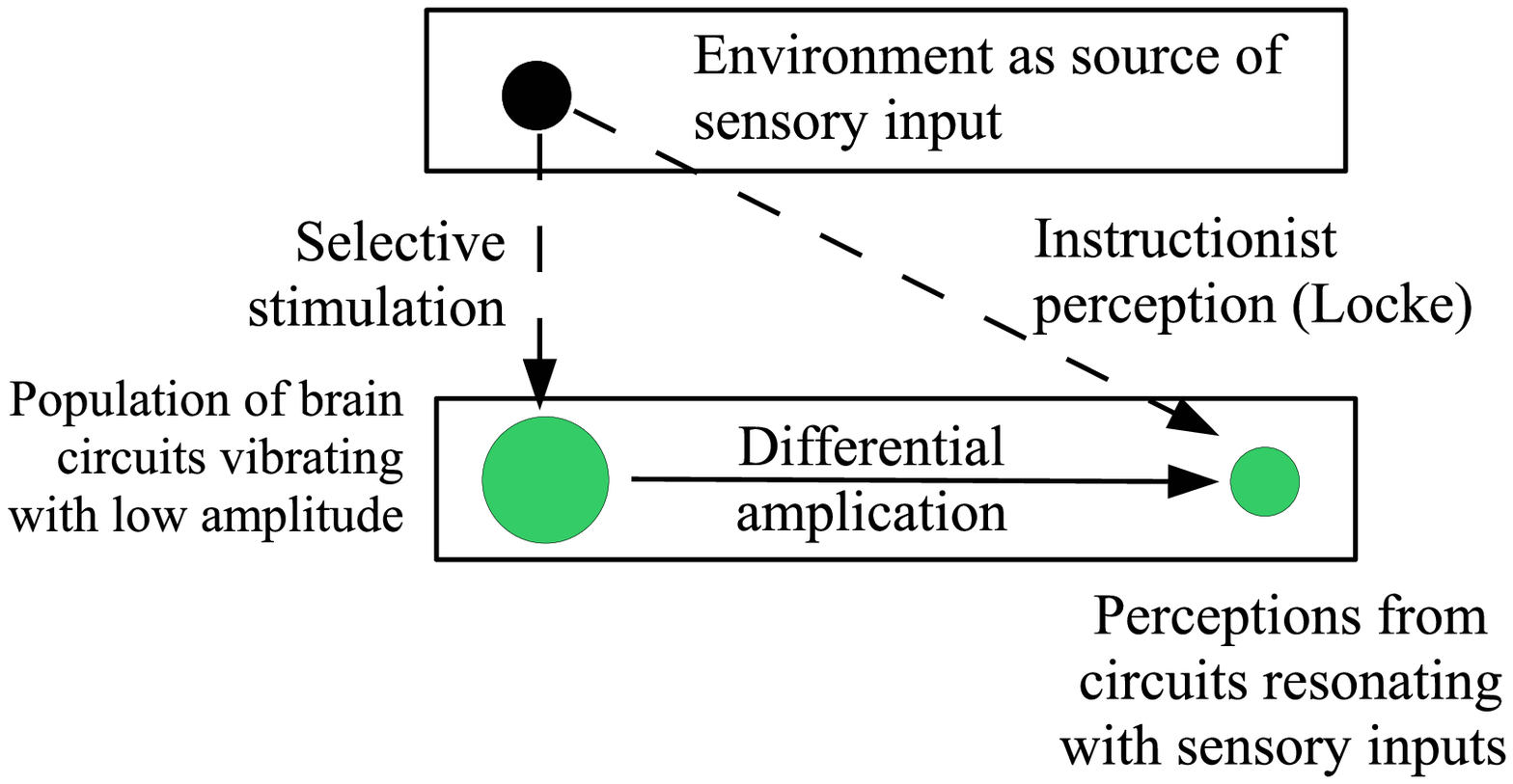}%
\\
Figure 8: Edelman's perception as determination through a receiving universal.
\end{center}

In Edelman's account of perception as the \textquotedblleft remembered
present,\textquotedblright\ direct determination is distinguished from the
composite effect of the indirect influence differentially triggering internal processes.

\begin{quote}
{\footnotesize According to this analysis, extrinsic signals convey
information not so much in themselves, but by virtue of how they modulate the
intrinsic signals exchanged within a previously experienced neural system. In
other words, a stimulus acts not so much by adding large amounts of extrinsic
information that need to be processed as it does by amplifying the intrinsic
information resulting from neural interactions selected and stabilized by
memory through previous encounters with the environment. \cite[p.
137]{edel:uc}}
\end{quote}

\subsection{Active versus passive learning}

This recalls a much older theme of "recollection." One of the tell-tale signs
of a process of determination through universals is the indirectness of the
factorization through a universal. Here again, an instructionist account might
be first given for a process that is later recognized as being selectionist.
The interplay between these two accounts dates back at least to the
Platonic-Socratic account of learning not as the result of external
instruction but as a process of catalyzing internal recollection. One of the
striking epigrams of neo-Platonism is the thesis that \textquotedblleft no man
ever does or can teach another anything\textquotedblright\ \cite[p.
1]{burn:witt}. In the early fifth century, Augustine in \textit{De Magistro
}(The Teacher) made the point contrasting \textquotedblleft
outward\textquotedblright\ passive instruction with active learning
\textquotedblleft within.\textquotedblright

\begin{quote}
{\footnotesize But men are mistaken, so that they call those teachers who are
not, merely because for the most part there is no delay between the time of
speaking and the time of cognition. And since after the speaker has reminded
them, the pupils quickly learn within, they think that they have been taught
outwardly by him who prompts them.(Chapter XIV)}
\end{quote}

\noindent In the nineteenth century, Wilheim von Humboldt made the same point
even recognizing the symmetry between listener and speaker (a symmetry
captured by the right semiadjunctions).

\begin{quote}
{\footnotesize Nothing can be present in the mind (Seele) that has not
originated from one's own activity. Moreover understanding and speaking are
but different effects of the selfsame power of speech. Speaking is never
comparable to the transmission of mere matter (Stoff). In the person
comprehending as well as in the speaker, the subject matter must be developed
by the individual's own innate power. What the listener receives is merely the
harmonious vocal stimulus. \cite[p. 102]{hum:herm}}
\end{quote}

\noindent A similar theme has been a mainstay in active learning theories of
education. As John Dewey put it:

\begin{quote}
{\footnotesize It is that no thought, no idea, can possibly be conveyed as an
idea from one person to another. When it is told, it is, to the one to whom it
is told, another given fact, not an idea. The communication may stimulate the
other person to realize the question for himself and to think out a like idea,
or it may smother his intellectual interest and suppress his dawning effort at
thought. \cite[p. 159]{dewey:de}}
\end{quote}

\noindent Remarkably, the immunologist Niels Jerne tied these themes together.

\begin{quote}
{\footnotesize Several philosophers, of course, have already addressed
themselves to this point. John Locke held that the brain was to be likened to
white paper, void of all characters, on which experience paints with almost
endless variety. This represents an instructive theory of learning, equivalent
to considering the cells of the immune system void of all characters, upon
which antigens paint with almost endless variety.}

{\footnotesize Contrary to this, the Greek Sophists, including Socrates, held
a selective theory of learning. Learning, they said, is clearly impossible.
For either a certain idea is already present in the brain, and then we have no
need of learning it, or the idea is not already present in the brain, and then
we cannot learn it either, for even if it should happen to enter from outside,
we could not recognize it. This argument is clearly analogous to the argument
for a selective mechanism for antibody formation, in that the immune system
could not recognize the antigen if the antibody were not already present.
Socrates concluded that all learning consists of being reminded of what is
pre-existing in the brain.\cite[pp. 204-5]{jerne:al}}
\end{quote}

\subsection{Chomsky's language acquisition faculty}

Niels Jerne also saw the connection with another major example of a left
semiadjunction, Norm Chomsky's theory of language acquisition. Language
learning by a child is another example of a process that was originally
thought to be instructive. \ But Noam Chomsky's theory of generative grammar
postulated an innate language faculty that would unfold according to the
linguistic experience of the child. \ The child did not \textquotedblleft
learn\textquotedblright\ the rules of grammar; the linguistic experience of
the child would select how the universal language faculty would develop or
unfold to differentially implement one rule rather than another. Niels Jerne
entitled his Nobel Lecture, \textit{The Generative Grammar of the Immune
System}.%

\begin{center}
\includegraphics[
height=2.3213in,
width=5.2304in
]%
{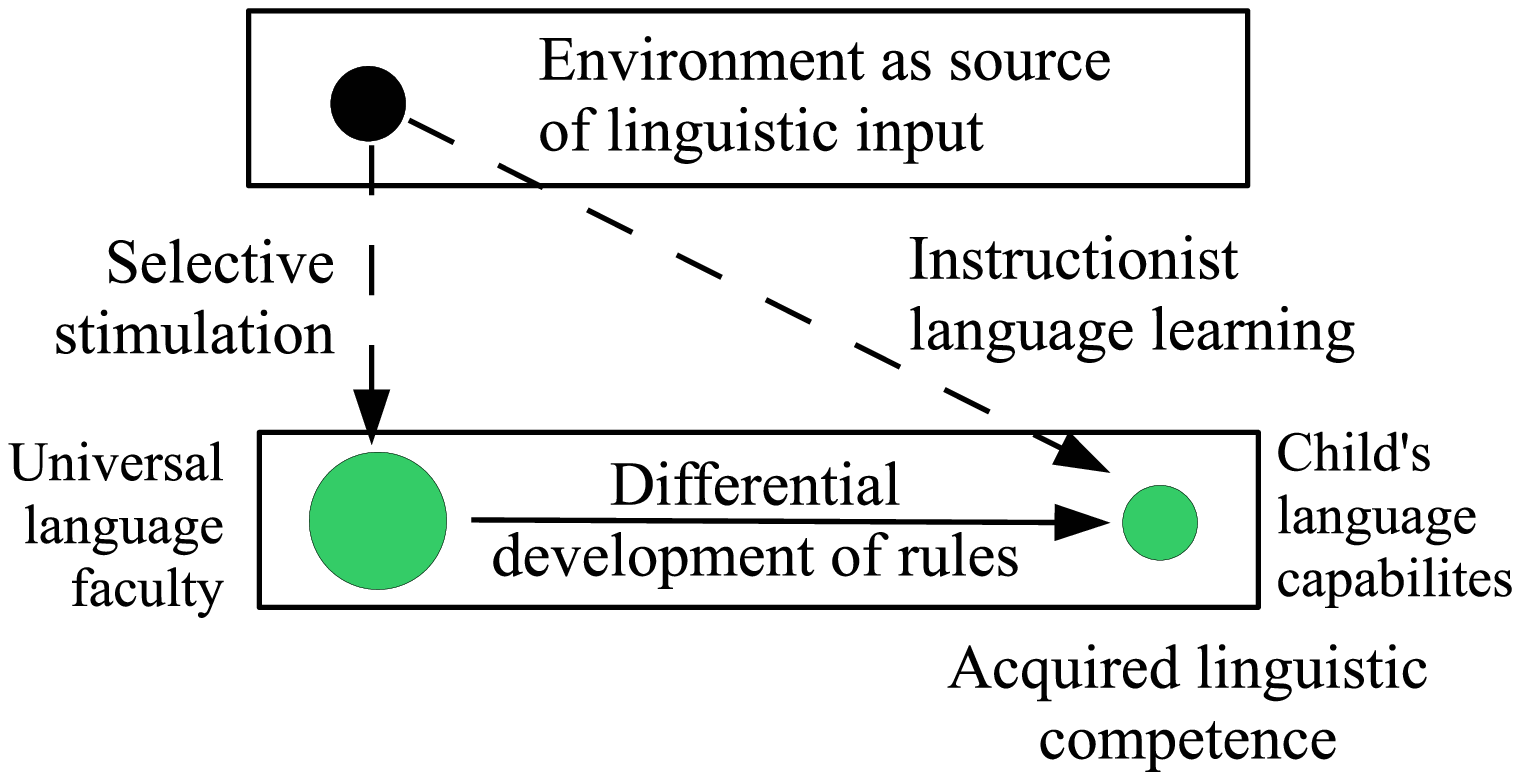}%
\\
Figure 9: Language acquisition as determination through a receiving universal.
\end{center}

Chomsky \cite{chom:cl} emphasizes certain aspects of the language acquisition
faculty which can be illustrated as general aspects of left semiadjunctions:

\begin{itemize}
\item \textit{Rich profligate internal structure}: the capacity to generate
`all' the various possibilities with the receiving universal;

\item \textit{Innateness}: that receiving universal is on the receiving
("organism") side, not the sending ("environment") side, even though it
"recognizes" the inputs from the sending side;

\item \textit{Impoverished or minimal inputs}: for instance, selecting the
jukebox button inputs less information that the whole instructive message of
the supplied record in Medawar's example;

\item \textit{Active role of internal mechanism}: the selection interacts with
the universal receiver to generate a specific internal determination rather
than passively receiving the detailed instructionist message like a stamp in wax;

\item \textit{Relative autonomy of internal mechanism}: the active internal
role together with the impoverished external input (and lack of
stimulus-control) gives the determination through the receiving universal a
certain type of autonomy from the sending environment.
\end{itemize}

Since one of the principal capabilities of Chomsky's language faculty is
recursion, we might make a slight mathematical digression to notice that the
mathematical formulation of recursion is a left semiadjunction.%

\begin{center}
\includegraphics[
height=1.7609in,
width=5.0909in
]%
{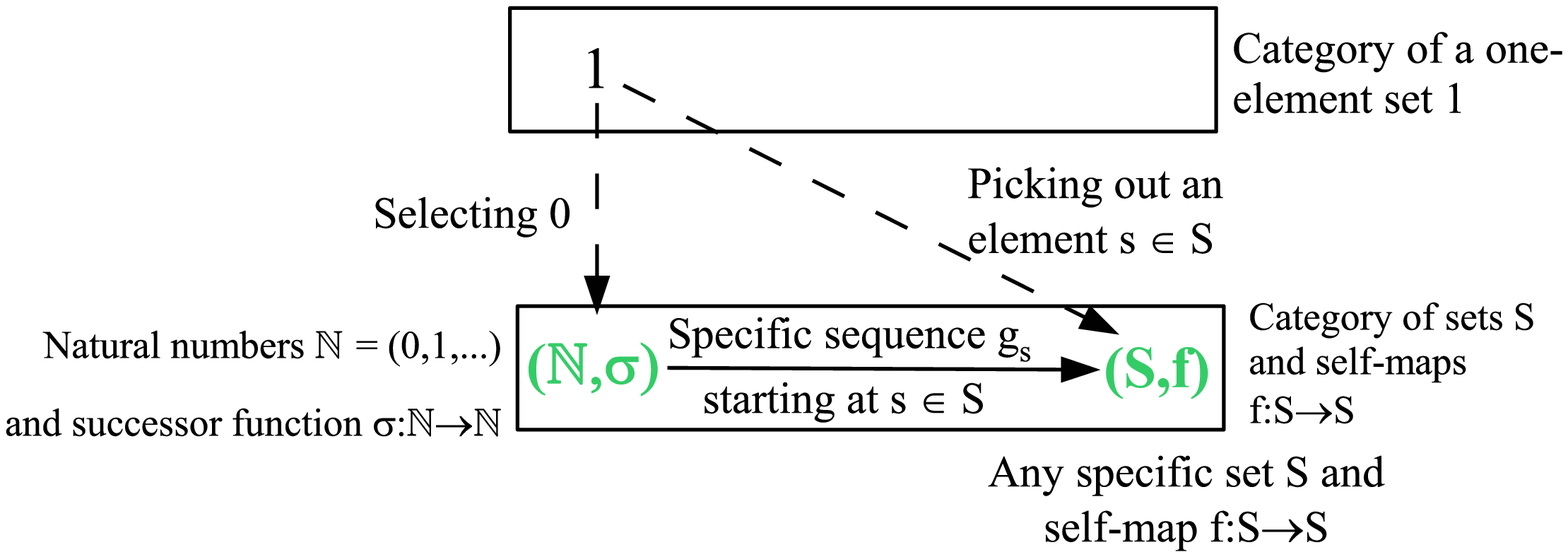}%
\\
Figure 10: Recursion as determination through a receiving universal.
\end{center}

\noindent The specific sequence $g_{s}$ enumerated by the natural numbers is
$s,f(s),f^{2}\left(  s\right)  ,f^{3}\left(  s\right)  ,...$.

\subsection{Language understanding}

Not only the acquisition of language but the ordinary understanding of
language is another example well-modeled by determination through a receiving
universal. If the auditory input is in a language that the listener
understands, that means there is an internal process triggered by the auditory
signals that recognizes, interprets, and understands the input. The Lockean or
behaviorist alternative is where the auditory input (het), Humboldt's "vocal
stimulus," is directly supplied by the auditory version of Locke's writing on
a blank slate.%

\begin{center}
\includegraphics[
height=2.4201in,
width=4.7887in
]%
{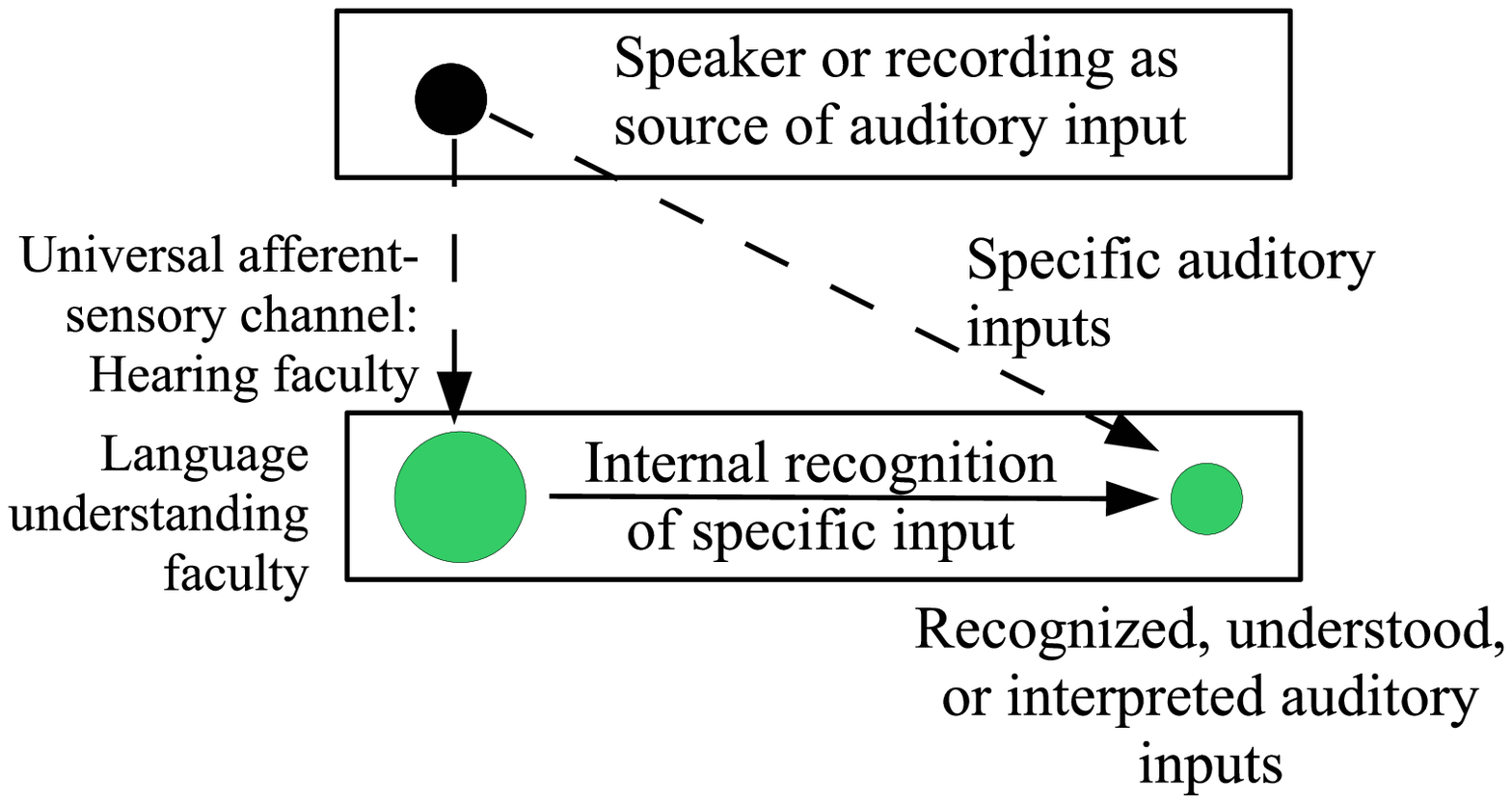}%
\\
Figure 11: Language understanding through a receiving universal.
\end{center}

\noindent The determination through the receiving universal thus adds a second
level to the input which is variously called recognition, understanding, or
the \textit{intentionality of perception}. This second level is often
indicated by saying that the sensory input is "perceived as", "recognized as",
or "understood as" something further.

\subsection{Generic "recognition" or "perception"}

Before turning to right semiadjunctions, it might be useful to present a
rather generic version of determination through a receiving universal as model
of "recognition" or "perception" that captures many of the common features of
the various examples.

The determination through the receiving universal is the active internal
process that supplies the "interpretation" or "intentionality" to the raw
sense data. The red blotch is seen as a tomato; the sound "ya" is understood
as indicating agreement, and so forth. In the direct alternative, the raw
sensory input supplies Lockean "perception" like writing on a blank slate or a
stamp making an impression on wax.%

\begin{center}
\includegraphics[
height=2.8427in,
width=4.9132in
]%
{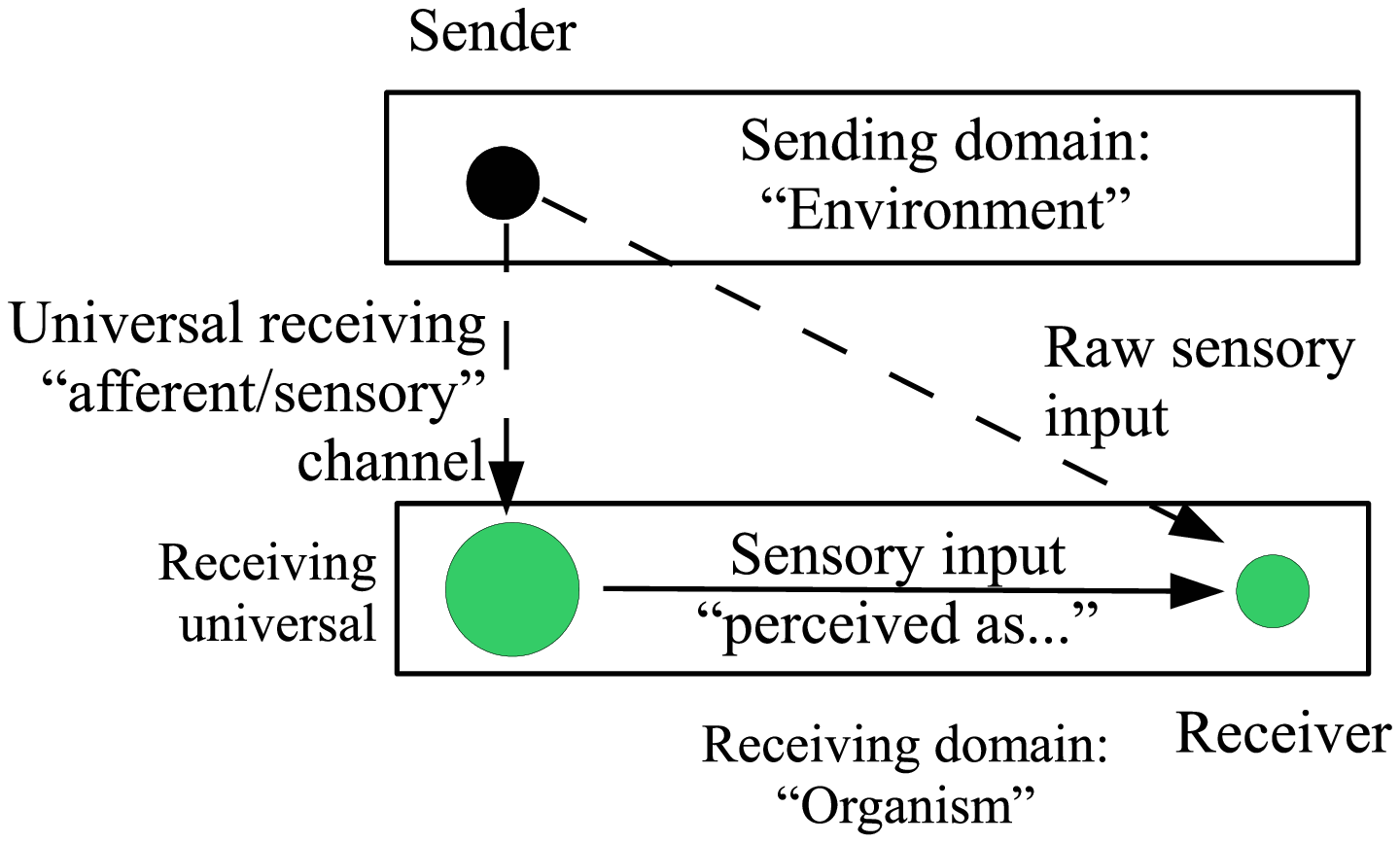}%
\\
Figure 12: "Recognition" as determination through a receiving universal.
\end{center}

\section{Determinations through a sending universal}

\subsection{Generic "action"}

Dual to the generic model of "recognition" is the generic model of
"action"--which is the determinative scheme given by a right semiadjunction.
In the model of recognition, there is the uninterpreted message as just a
sensory input (the external het), and then there is the second level where the
factorization (the internal hom) through the receiving universal recognizes
the interpretation, meaning, or intentionality of the message. In the dual
model of "action," the external het specifies the external behavior (which
could be even a reflex behavior) while internal hom factoring through a
sending universal that supplies the "intentionality" of the "action" (where an
"action" is a "behavior" with the second level of "intentionality"). In each
case, we end up with a certain behavior but determined by two different means.

\begin{center}%
\begin{center}
\includegraphics[
height=2.0398in,
width=4.4549in
]%
{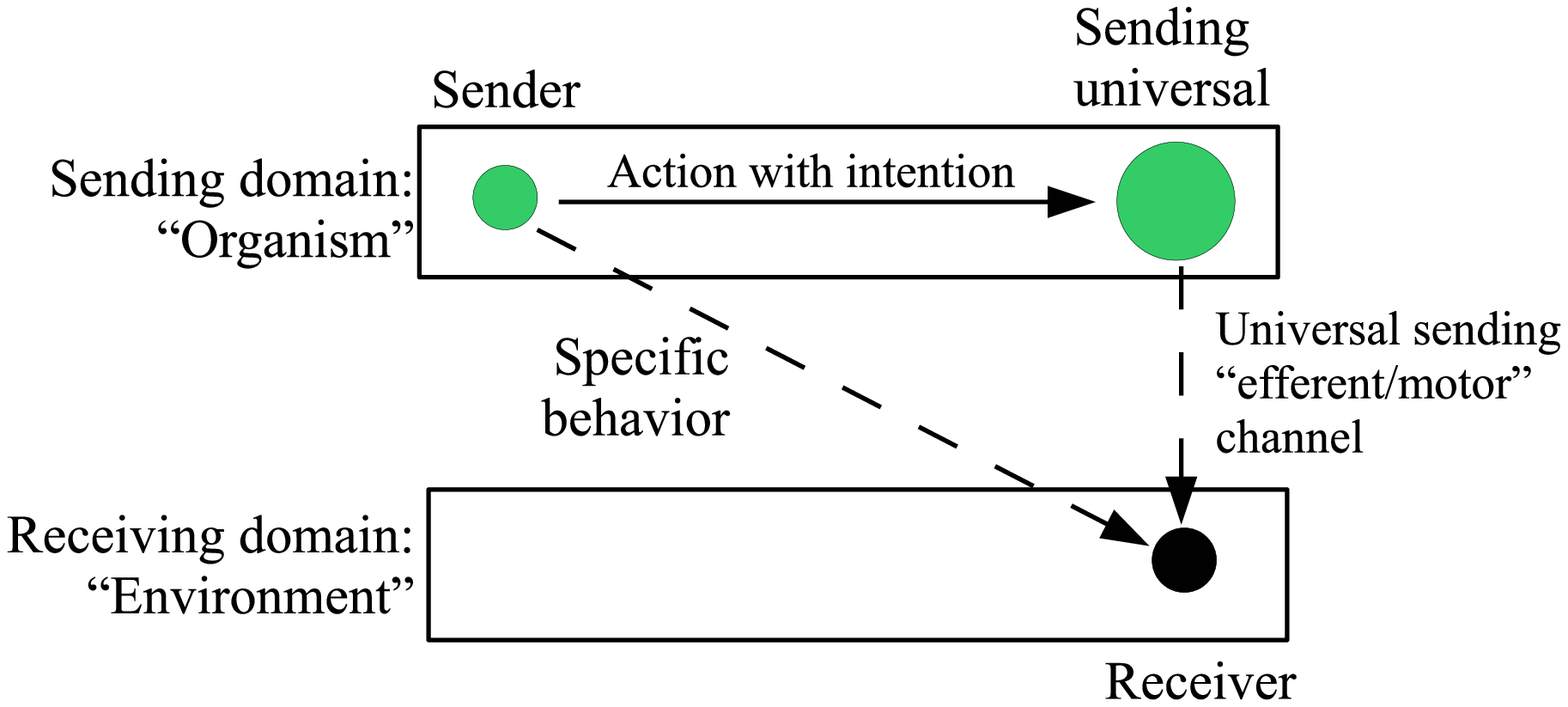}%
\\
Figure 13: "Action" as determination through a sending universal.
\end{center}

\end{center}

\subsection{Universal Turing machines}

Universal Turing machines provide a good example of a sending universal and
they are strikingly similar to the next example of the genetic code. Given the
inputs, a specific Turing machine (TM) might calculate a function such as the
successor function. But if those inputs and the coding for that specific TM
were fed in as inputs to a universal TM, then it would realize the same
calculated results.%

\begin{center}
\includegraphics[
height=2.7738in,
width=4.2881in
]%
{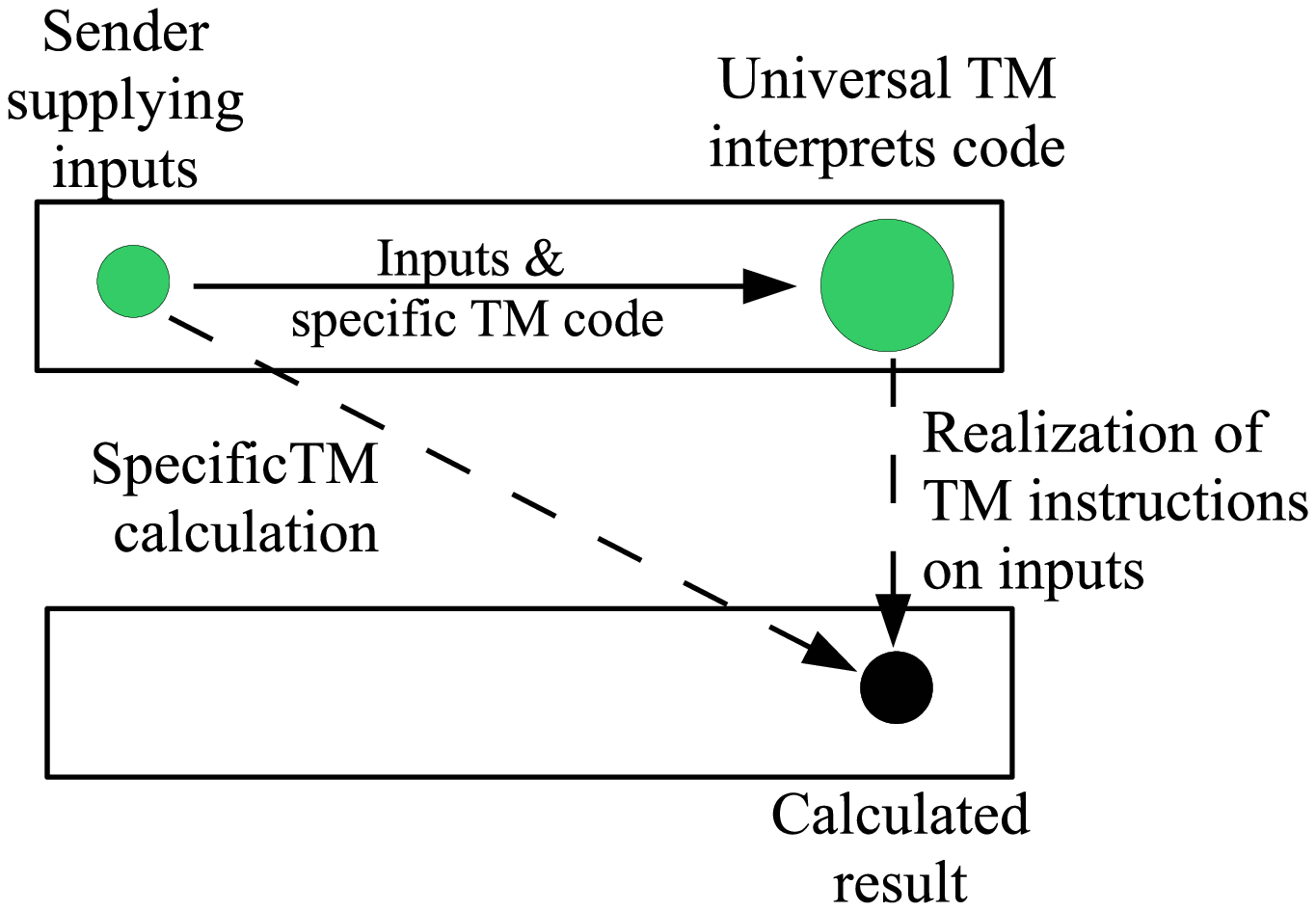}%
\\
Figure 14: Universal Turing machine as a sending universal.
\end{center}

\subsection{The genetic code and DNA mechanism}

Surely the most important example of a right semiadjunction in biology is the
whole DNA mechanism to use the genetic code to construct amino acids. One
might imagine a specific chemical mechanism that would produce a specific
amino acid--which would play the role of the specific sending het. But the DNA
mechanism uses an internal sending universal to interpret coded messages to
make one amino acid or another, and then to realize the construction of that
amino acid.%

\begin{center}
\includegraphics[
height=3.1623in,
width=5.0211in
]%
{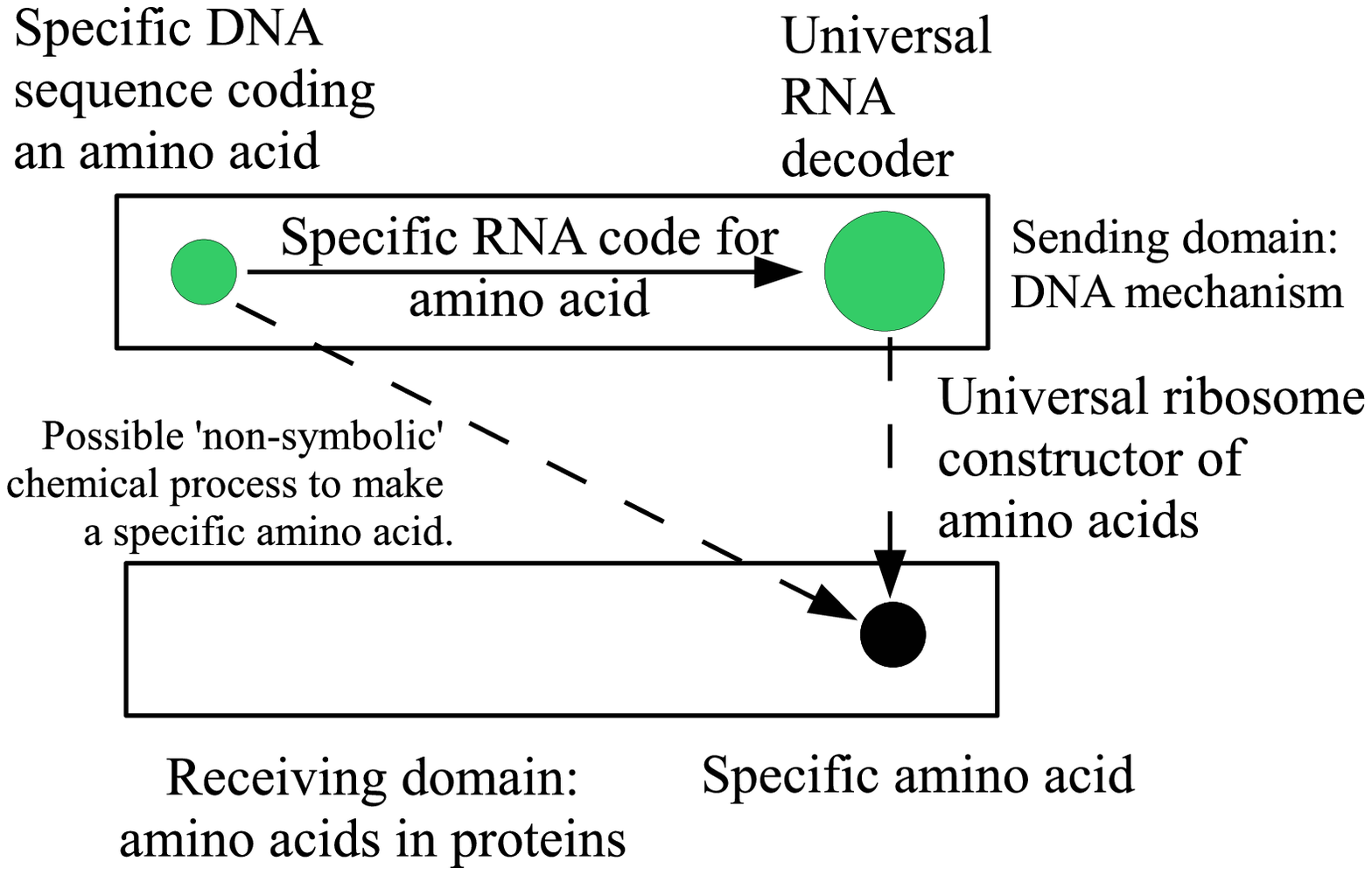}%
\\
Figure 15: DNA mechanism as a determination through a sending universal.
\end{center}

John Maynard Smith has emphasized that "the use of informational terms implies
intentionality" \cite[p. 124]{maynard:info} which is that second level, the
chemistry plus the coded meaning of the chemistry. The "meaning" is supplied
not by a human sender but by the process of evolution which essentially sends
the message: "this form has survival value."

Furthermore it is the use of coded information that allows for the generation
of such "universal" variety for selection to operate upon:

\begin{quote}
I think that it is the symbolic nature of molecular biology that makes
possible an indefinitely large number of biological forms. \cite[p.
133]{maynard:info}
\end{quote}

\begin{quotation}
The mere existence of replicators is not sufficient for continued evolution by
natural selection, which requires what we have called "indefinite hereditary
replicators": that is, entities that can exist in an indefinitely large number
of states, each of which can be replicated. In living organisms, nucleic acid
molecules are the only indefinite hereditary replicators, or at least they
were until the invention of language and music. \cite[p. 58]%
{maynard;transitions}
\end{quotation}

\subsection{Developmental mechanisms}

Less well understood is the way that informational codes drive hierarchies of
regulator genes to develop this or that organ. But the whole arrangement still
has the form of a right semiadjunction. The contrast between a specific het
and the corresponding specific hom (with its second level of meaning or
intentionality) is explained by Maynard Smith using the example of proteins.

\begin{quotation}
A protein might have a function directly determined by its structure--for
example, it may be a specific enzyme, or a contractile fiber. Alternatively,
it might have a regulatory function, switching other genes on or off. Such
regulatory functions are arbitrary, or symbolic. They depend on specific
receptor DNA sequences, which have themselves evolved by natural selection.
The activity of an enzyme depends on the laws of chemistry and on the chemical
environment (for example, the presence of a suitable substrate), but there is
no structure that can be thought of as an evolved "receiver" of a "message"
from the enzyme. By contrast, the effect of a regulatory protein does depend
on an evolved receiver of the information it carries: the eyeless gene signals
"make an eye here," but only because the genes concerned with making an eye
have an appropriate receptor sequence. \cite[p. 143]{maynard:info}
\end{quotation}

%

\begin{center}
\includegraphics[
height=2.5712in,
width=4.3288in
]%
{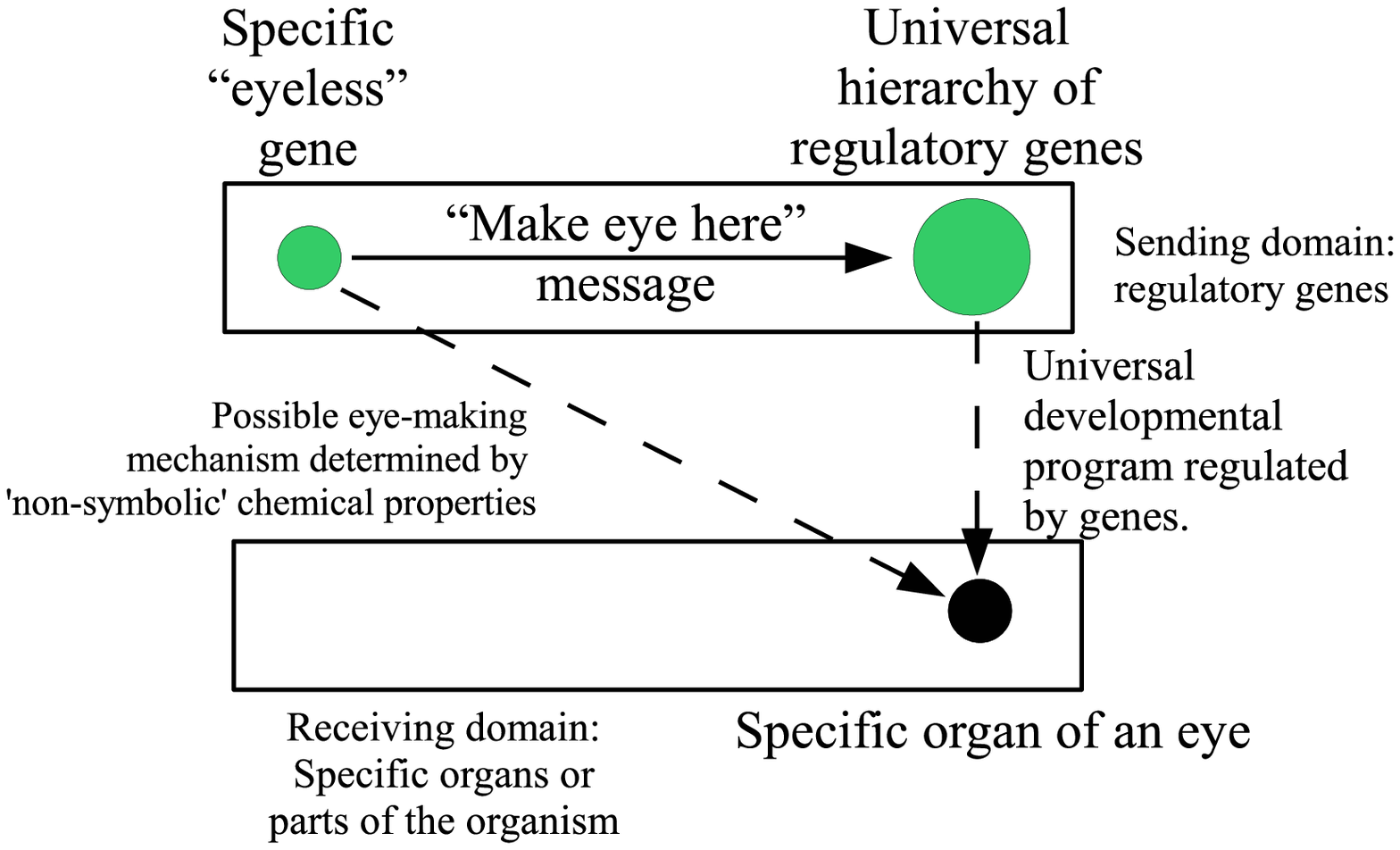}%
\\
Figure 16: Hierarchy of regulatory genes as a sending universal.
\end{center}

A stem cell would fit into the same type of diagram as the sending universal.
One could imagine some specific non-symbolic chemical process (playing the
role of the specific het) that would develop a certain type of bodily cell.
But a stem cell is a sending universal in that upon receiving the proper
specific chemical codes (playing the role of the corresponding specific hom),
it will develop into some specific type of cell.

\subsection{Language action}

The dual to "language understanding" is language production or linguistic
action (e.g., "speech acts"). The role of the specific het is played by some
auditory output such as utterances (Humboldt's "vocal stimulus"). But the
corresponding internal specific hom is the speech act (i.e., internal speech
with intentionality) that through the language faculty produces the same
outputs but as intentional speech.%

\begin{center}
\includegraphics[
height=2.1818in,
width=4.0573in
]%
{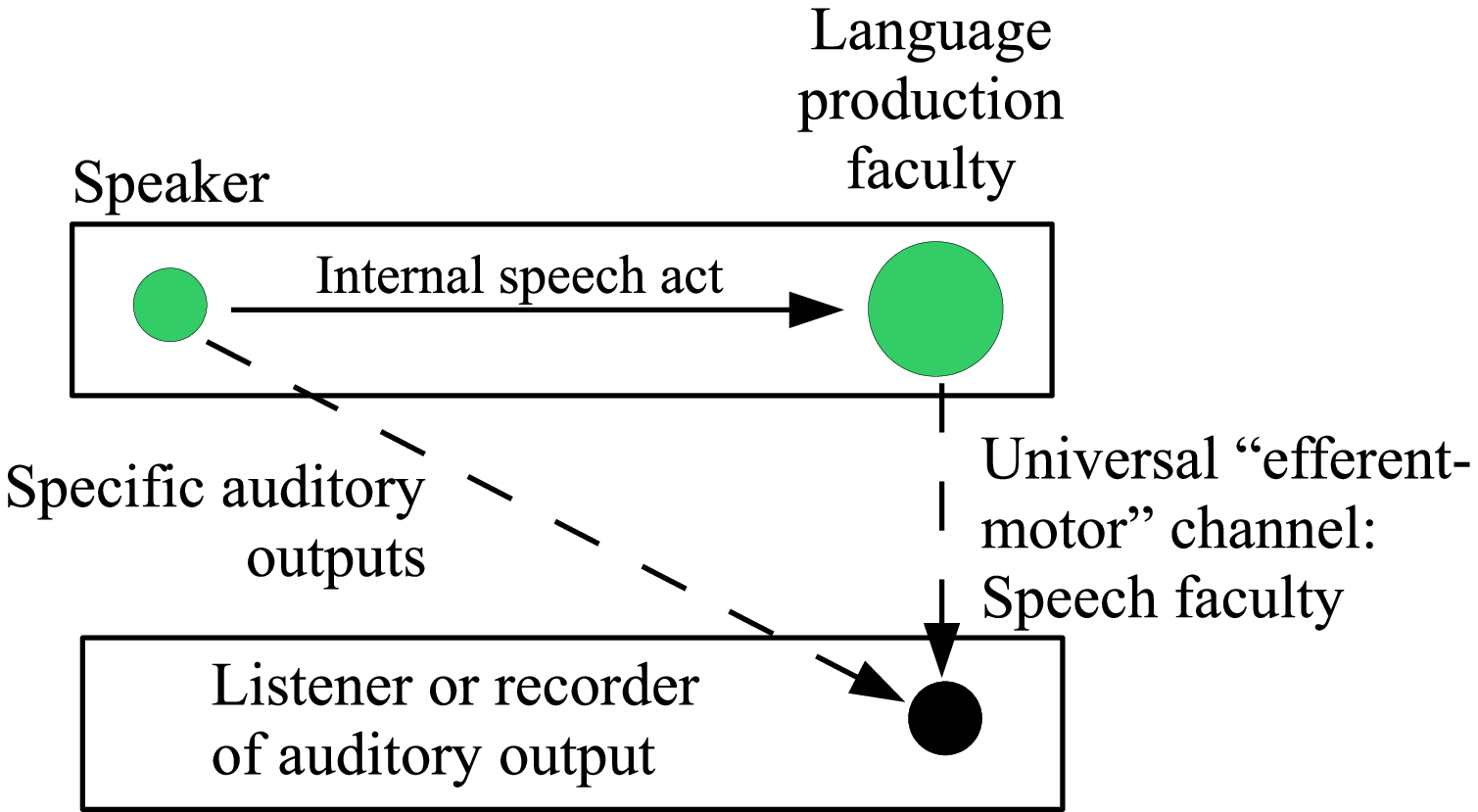}%
\\
Figure 17: Language production through a sending universal.
\end{center}

\section{Brain Functors}

\subsection{Recombination of left and right semiadjunctions}

It might be noted that in a number of our examples essentially the same
faculty, e.g., the language faculty, plays both the role of the receiving or
recognition universal and the sending or action universal. Hence it should be
possible to recombine the semiadjunctions so that the two universals coincide,
and that yields the concept of the brain functor (the "brain" being that
universal for both recognition and action). From the mathematical viewpoint,
the key to this was using the hets to split an adjunction into left and right
semiadjunctions which can then be recombined in the opposite way.%

\begin{center}
\includegraphics[
height=3.5085in,
width=4.7679in
]%
{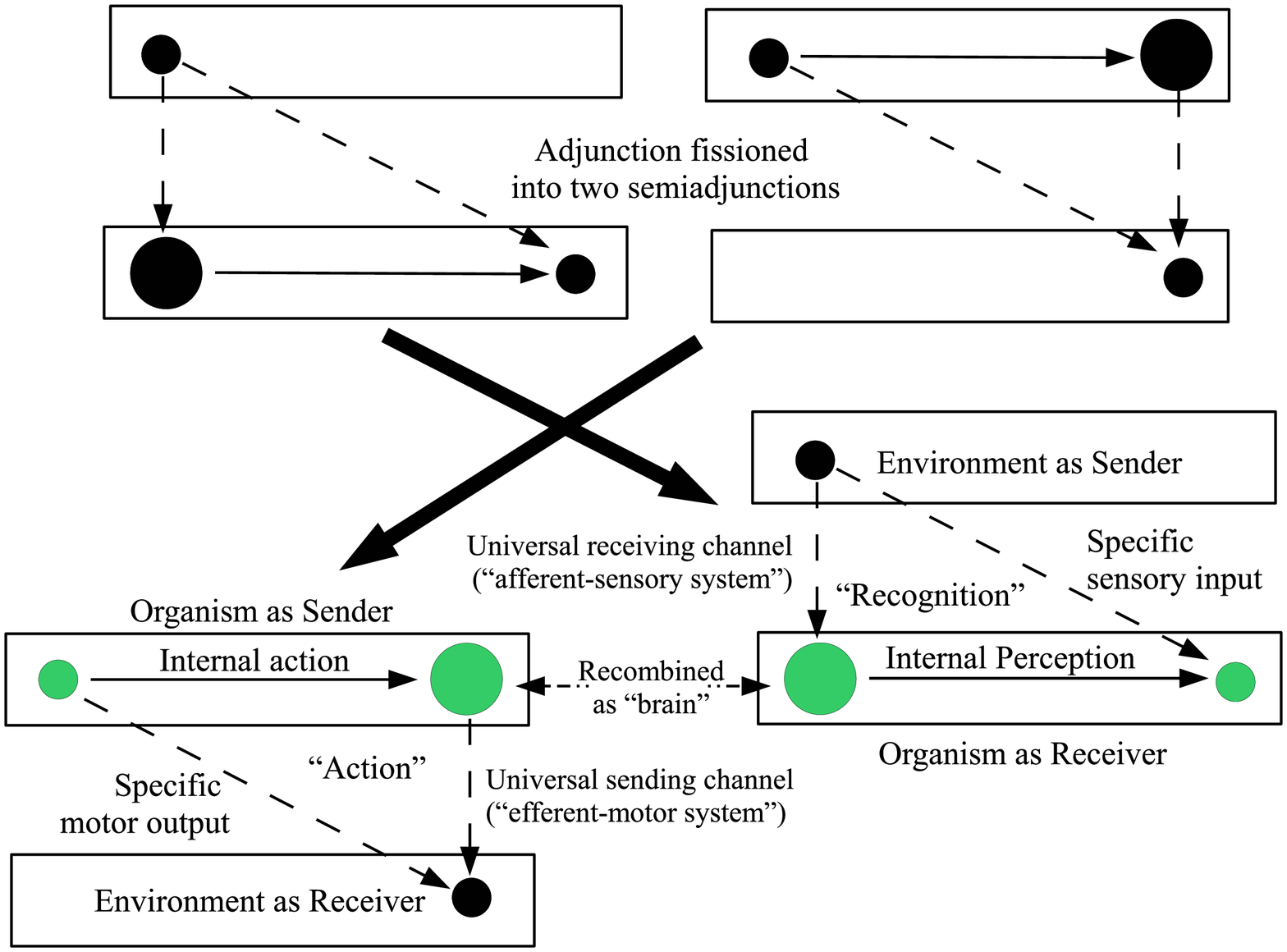}%
\\
Figure 18: Recombining a fissioned adjunction to make a brain functor.
\end{center}

The recombined semiadjunctions then form a \textit{brain functor}. In the
following `butterfly' for a brain functor, we use labels appropriate to the
"brain" label.%

\begin{center}
\includegraphics[
height=3.276in,
width=5.0087in
]%
{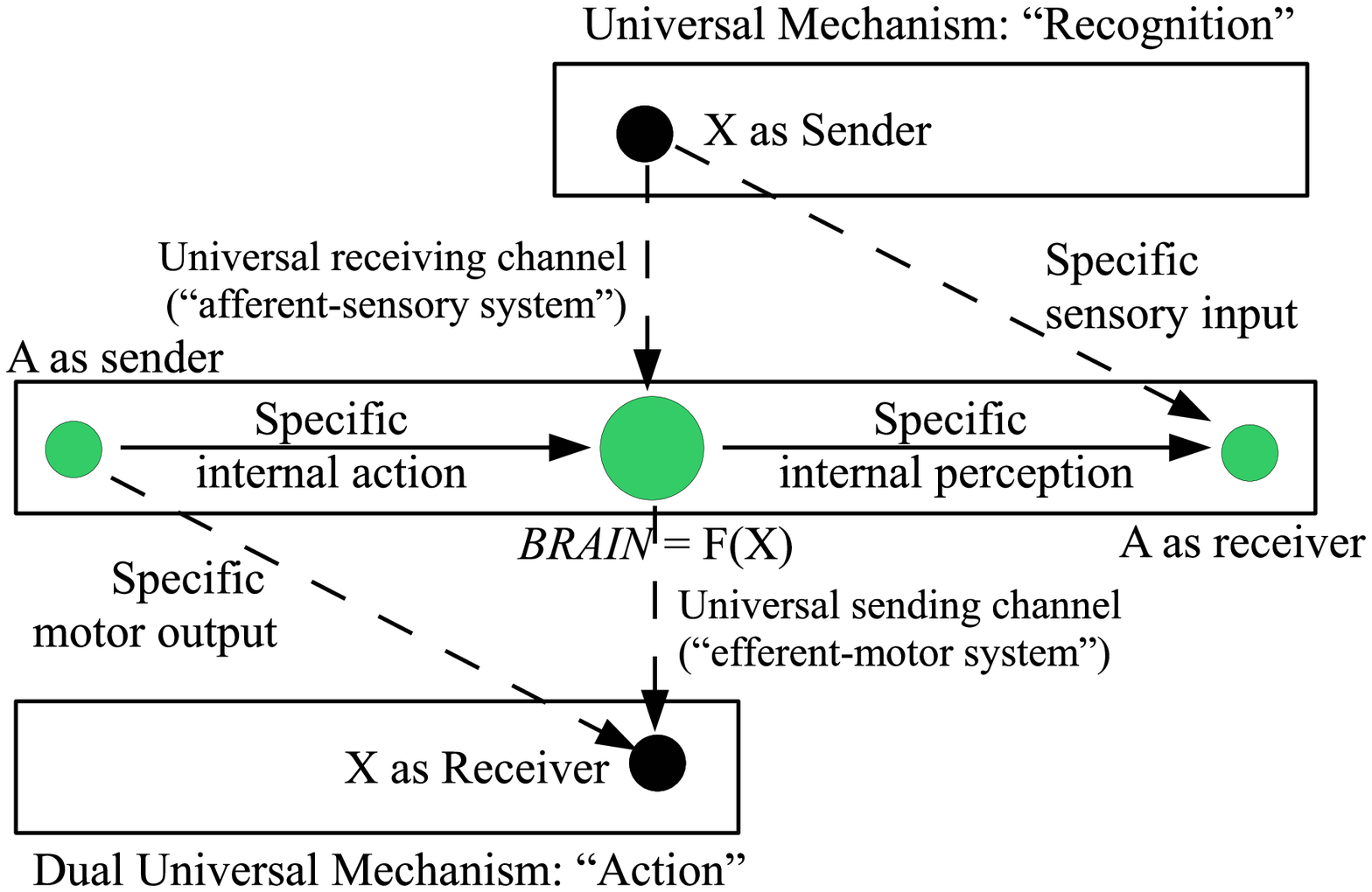}%
\\
Figure 19: Brain: Scheme for receiving and sending through one universal.
\end{center}

Mathematically, for each het $X\dashrightarrow A$ there is a unique internal
"recognizing" hom $F(X)\rightarrow A$ so there is a canonical isomorphism:
$Hom(F(X),A)\cong Het(X,A)$. And for each het $A\dashrightarrow X$ in the
other direction, there is a unique internal "action" hom $A\rightarrow
F\left(  X\right)  $ so there is also a canonical isomorphism: $Het(A,X)\cong
Hom(A,F(X))$.

The concept of a brain functor is the natural cognate or associated concept to
an adjunction. For an adjunction, there are \textit{two} functors that
represent on the left and right the hets going \textit{one} way between the categories:

\begin{center}
$Hom(F(X),A)\cong Het(X,A)\cong Hom(X,G(A))$.
\end{center}

\noindent For a brain functor, there is \textit{one} functor that represents
on the left and right the hets going the \textit{two} ways between the categories:

\begin{center}
$Hom(F\left(  X\right)  ,A)\cong Het(X,A)$ and $Het(A,X)\cong Hom(A,F(X))$.
\end{center}

\noindent Hence the general scheme given by a brain functor is
\textit{receiving and sending determination through one universal} (the "brain").

\subsection{Simple two-way determinations through one universal}

The simplest form of a brain "functor" is just a two-way representation or
coding system that constructs and implements a set of codes. Given some set of
objects, it is encoded using some isomorphic set of representations or codes
for the objects, and then given an instance of the code, it is decoded to
determine the object.

Coordinatizing is a form of coding. The geometrical plane is a collection of
points, and the Cartesian coordinate system represents each point $P$ by a
pair $\left(  x_{P},y_{P}\right)  $ of coordinates. Given a point $P$, the
"coordinate" function selects the coordinates $\left(  x_{P},y_{P}\right)  $
of the point which is the recognized or coded output, and given the
coordinates or code for a point $\left(  x_{P},y_{P}\right)  $ as an input,
the "plot" function designates the point.%

\begin{center}
\includegraphics[
height=2.5529in,
width=5.271in
]%
{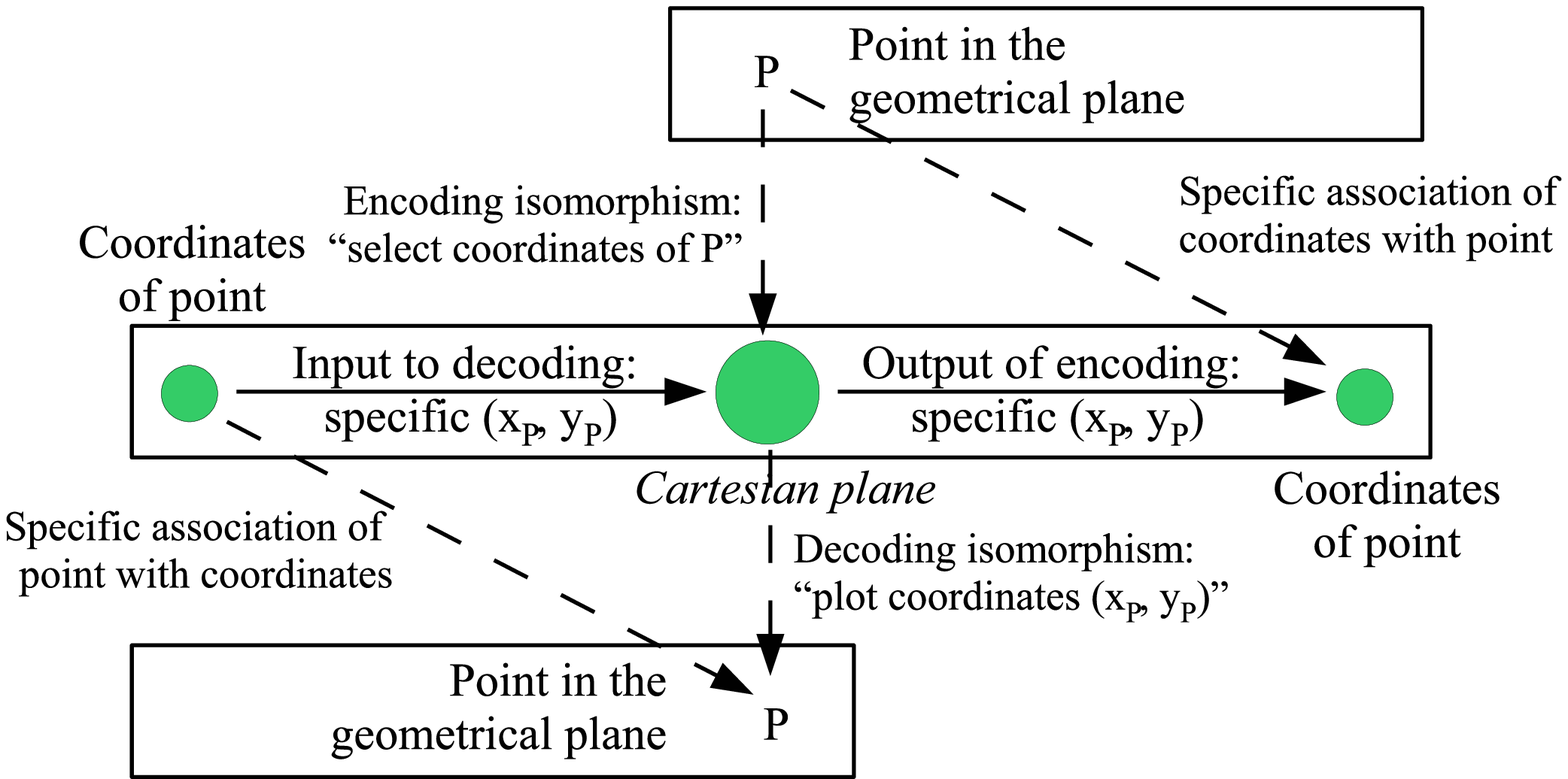}%
\\
Figure 20: Coding and decoding Cartesian coordinates of geometrical points.
\end{center}

\subsection{The language faculty}

A brain functor, broadly put, is any universal mechanism of determination that
can factor determination either way through a universal--rather than an
adjunction that factors one way determination through two (receiving and
sending) universals. In some contexts in the life sciences, determination is
strictly one way so one might expect to find a semiadjunction but not a
two-way system like a brain functor. For instance, there is the
\textit{fundamental dogma} that DNA determines amino acids, proteins, and
eventually the characteristics of an organism but never the reverse.

An application of the scheme for a brain functor in the cognitive sciences is
to model the language faculty where there is two way determination between
vocal stimuli and internal representations. The previous semiadjunctions for
language understanding and language action can be merged to arrive at the
brain-like function of the language faculty.%

\begin{center}
\includegraphics[
height=2.562in,
width=5.139in
]%
{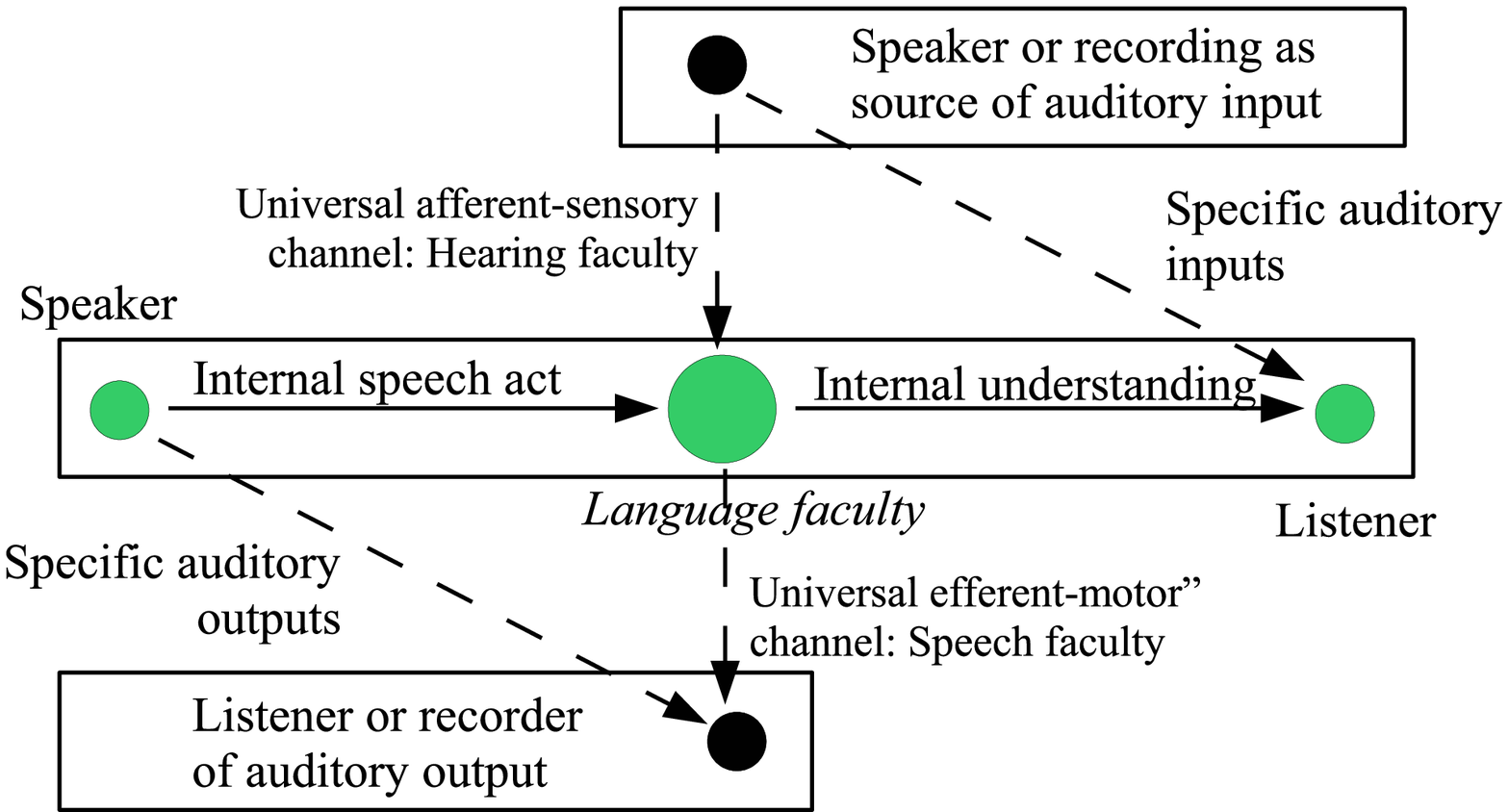}%
\\
Figure 21: Language faculty as two-way determination through a universal.
\end{center}

\section{Summary}

The following table gives the principal exact category-theoretic concepts to
describe the corresponding schemes of determination through universals as well
as the main generic examples.

\begin{center}%
\begin{tabular}
[c]{|c|c|c|}\hline
\textbf{CT concept} & \textbf{Determination through universals} &
\textbf{Generic example}\\\hline\hline
Left semiadjunction & Det. through a receiving universal & Recognition\\\hline
Right semiadjunction & Det. through a sending universal & Action\\\hline
Brain functor & Two-way det. through a universal & Recognition +
Action\\\hline
\end{tabular}

Table 1: Principal forms of determination through universals
\end{center}

The importance of category-theoretical universals in isolating the important
concepts in pure mathematics suggests that the universals may play a similar
role in the empirical sciences. Our results suggest that this is indeed the
case in the biological and cognitive sciences. The category-theoretic schemes
of determination through universals are at a high level of abstraction, but,
in this case at least, that seems to be where some significant theory lives.
Regardless of the great differences in the underlying substrate processes,
many of the most important mechanisms and faculties of the biological and
cognitive sciences not only seem to fit into, but also to have their key
features characterized by, the schemes for determination through universals.

\section{Appendix: Defining hets in category theory}

Category theory groups together in \textit{categories} the mathematical
objects with some common structure (e.g., sets, partially ordered sets,
groups, rings, and so forth) and the appropriate morphisms between such
objects. Since the morphisms are between objects of similar structure, they
are ordinarily called \textquotedblleft homomorphisms.\textquotedblright\ 

But there have always been other morphisms which occur in mathematical
practice that are between objects with different structures (i.e., in
different categories) such as the insertion-of-generators map from a set to
the free group on that set. Indeed, the working mathematician might well
characterize the free group $F\left(  X\right)  $ on a set $X$ as the group
such that for any set-to-group map $f:X\dashrightarrow G$, there is a unique
group homomorphism $\widehat{f}:F\left(  X\right)  \rightarrow G$ that factors
$f$ \ through the canonical insertion of generators $i:X\dashrightarrow
F\left(  X\right)  $, i.e., $f=\widehat{f}i$. In order to contrast these
morphisms such as $f:X\dashrightarrow G$ and $i:X\dashrightarrow F\left(
X\right)  $ with the homomorphisms between objects within a category such as
$\widehat{f}:F\left(  X\right)  \rightarrow G$, the former are called
\textit{heteromorphisms} or \textit{hets} (for short). Hets are like chimeras
since they have a tail in one category and a head in another category.

We assume familiarity with the usual machinery of category theory (bifunctors,
in particular) which can be adapted to give a rigorous treatment of
heteromorphisms (and their compositions with homomorphisms) that is parallel
to the usual bifunctorial treatment of homomorphisms.

The cross-category object-to-object hets $d:X\dashrightarrow A$ will be
indicated by dashed arrows ($\dashrightarrow$) rather than solid arrows
($\rightarrow$). The first question is how do heteromorphisms compose with one
another? But that is not necessary. Chimera do not need to `mate' with other
chimera to form a `species' or category; they only need to mate with the
intra-category morphisms on each side to form other chimera.\footnote{The
chimera genes are dominant in these mongrel matings. While mules cannot mate
with mules, it is `as if' mules could mate with either horses or donkeys to
produce other mules.}

Given a het $d:X\dashrightarrow A$ from an object in a category $\mathbf{X}$
to an object in a category $\mathbf{A}$, and homs $h:X^{\prime}\rightarrow X$
in $\mathbf{X}$ and $k:A\rightarrow A^{\prime}$ in $\mathbf{A}$, the
composition $dh:X^{\prime}\rightarrow X\dashrightarrow A$ is another het
$X^{\prime}\dashrightarrow A$ and the composition $kd:X\dashrightarrow
A\rightarrow A^{\prime}$ is another het $X\dashrightarrow A^{\prime}$.%

\begin{center}
\includegraphics[
height=1.7609in,
width=4.5496in
]%
{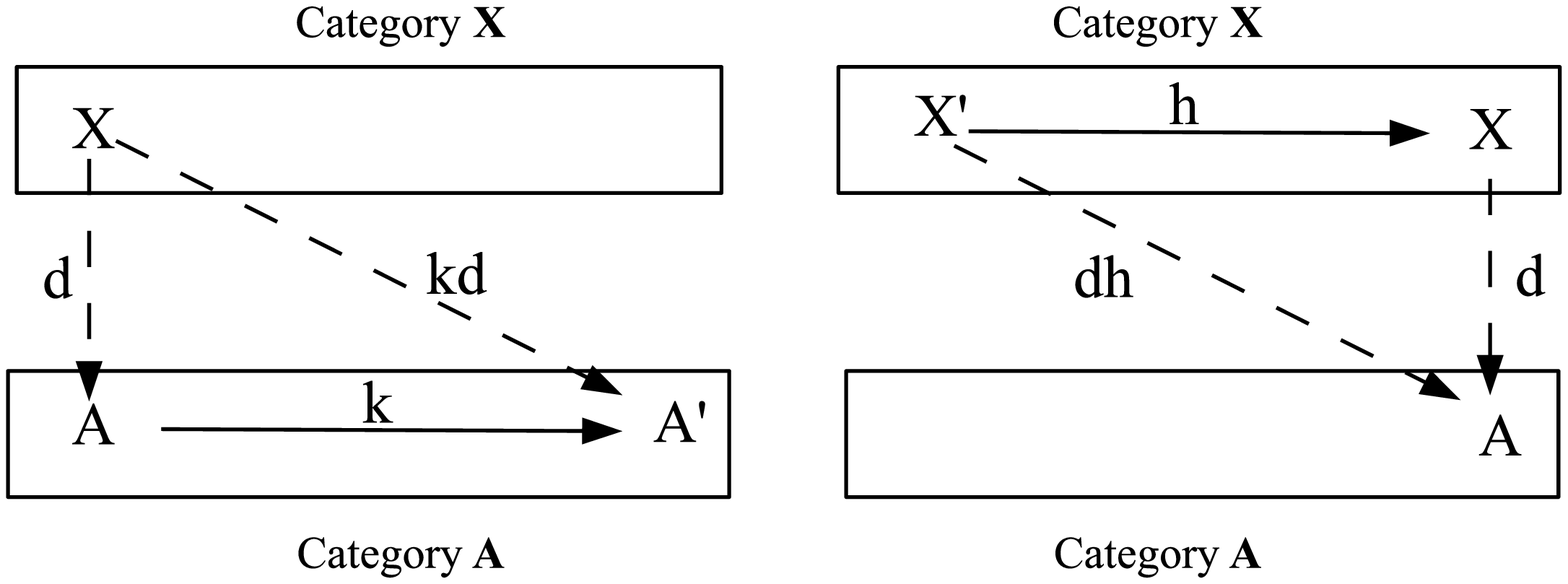}%
\\
Figure 22: Composition of hets and homs.
\end{center}

This action is exactly described by a bifunctor $\operatorname*{Het}%
:\mathbf{X}^{op}\times\mathbf{A}\rightarrow\mathbf{Set}$ where
$\operatorname*{Het}(X,A)=\{X\dashrightarrow A\}$ and where $\mathbf{Set}$ is
the category of sets and set functions. The natural machinery to treat
object-to-object morphisms \textit{between} categories are het-bifunctors
$\operatorname*{Het}:\mathbf{X}^{op}\times\mathbf{A}\rightarrow\mathbf{Set}$
that generalize the hom-bifunctors $\operatorname*{Hom}:\mathbf{X}^{op}%
\times\mathbf{X}\rightarrow\mathbf{Set}$ used to treat object-to-object
morphisms \textit{within} a category.

For any $\mathbf{A}$-hom $k:A\rightarrow A^{\prime}$ and any het
$X\overset{d}{\dashrightarrow}A$ in $\operatorname*{Het}(X,A)$, there is a
composite het $X\overset{d}{\dashrightarrow}A\overset{k}{\rightarrow}%
A^{\prime}=X\overset{kd}{\dashrightarrow}A^{\prime}$, i.e., $k$ induces a map
$\operatorname*{Het}(X,k):\operatorname*{Het}(X,A)\rightarrow
\operatorname*{Het}(X,A^{\prime})$. For any $\mathbf{X}$-hom $h:X^{\prime
}\rightarrow X$ and het $X\overset{d}{\dashrightarrow}A$ in
$\operatorname*{Het}(X,A)$, there is the composite het $X^{\prime
}\overset{h}{\rightarrow}X\overset{d}{\dashrightarrow}A=X^{\prime
}\overset{dh}{\dashrightarrow}A$, i.e., $h$ induces a map $\operatorname*{Het}%
(h,A):\operatorname*{Het}(X,A)\rightarrow\operatorname*{Het}(X^{\prime},A)$
(note the reversal of direction). The induced maps would respect identity and
composite morphisms in each category. Moreover, composition is associative in
the sense that $(kd)h=k(dh)$. This means that the assignments of sets of
chimera morphisms $\operatorname*{Het}(X,A)=\{X\overset{d}{\dashrightarrow
}A\}$ and the induced maps between them constitute a \textit{bifunctor}
$\operatorname*{Het}:\mathbf{X}^{op}\times\mathbf{A}\rightarrow\mathbf{Set}$
(contravariant in the first variable and covariant in the second).

With this motivation, we may turn around and define \textit{heteromorphisms}%
\ from $\mathbf{X}$-objects to $\mathbf{A}$-objects as the elements in the
values of a given bifunctor $\operatorname*{Het}:\mathbf{X}^{op}%
\times\mathbf{A}\rightarrow\mathbf{Set}$. This would be analogous to defining
the homomorphisms in $\mathbf{X}$ as the elements in the values of a given
hom-bifunctor $\operatorname*{Hom}_{\mathbf{X}}:\mathbf{X}^{op}\times
\mathbf{X}\rightarrow\mathbf{Set}$ and similarly for $\operatorname*{Hom}%
_{\mathbf{A}}:\mathbf{A}^{op}\times\mathbf{A}\rightarrow\mathbf{Set}$.

Given any bifunctor $\operatorname*{Het}:\mathbf{X}^{op}\times\mathbf{A}%
\rightarrow\mathbf{Set}$, it is \textit{representable on the left} if for each
$\mathbf{X}$-object $X$, there is an $\mathbf{A}$-object $F(X)$ that
represents the functor $\operatorname*{Het}(X,-)$, i.e., there is an
isomorphism $\psi_{X,A}:\operatorname*{Hom}_{\mathbf{A}}(F(X),A)\cong%
\operatorname*{Het}(X,A)$ natural in $A$. This defines a functor
$F:\mathbf{X}\rightarrow\mathbf{A}$, and any such functor with natural
isomorphisms $\operatorname*{Hom}_{\mathbf{A}}(F(X),A)\cong\operatorname*{Het}%
(X,A)$ is a \textit{left semiadjunction}.

Given a bifunctor $\operatorname*{Het}:\mathbf{X}^{op}\times\mathbf{A}%
\rightarrow\mathbf{Set}$, it is \textit{representable on the right} if for
each $\mathbf{A}$-object $A$, there is an $\mathbf{X}$-object $G(A)$ that
represents the functor $\operatorname*{Het}(-,A)$, i.e., there is an
isomorphism $\varphi_{X,A}:\operatorname*{Het}(X,A)\cong\operatorname*{Hom}%
_{\mathbf{X}}(X,G(A))$ natural in $X$. This defines a functor $G:\mathbf{A}%
\rightarrow\mathbf{X}$, and any such functor with natural isomorphisms
$\operatorname*{Het}(X,A)\cong\operatorname*{Hom}_{\mathbf{X}}(X,G(A))$ is a
\textit{right semiadjunction}.

An \textit{adjunction} is given by two functors $F:\mathbf{X}\rightarrow
\mathbf{A}$ and $G:\mathbf{A}\rightarrow\mathbf{X}$ that form left and right
semiadjunctions, i.e.,

\begin{center}
$\operatorname*{Hom}_{\mathbf{A}}(F(X),A)\cong\operatorname*{Het}%
(X,A)\cong\operatorname*{Hom}_{\mathbf{X}}(X,G(A))$
\end{center}

\noindent natural for all $X,A$. The identity hom $1_{F(X)}\in
\operatorname*{Hom}_{\mathbf{A}}(F(X),F(X))$ is associated by the left
isomorphism with the \textit{universal receiving het} $h_{X}:X\dashrightarrow
F(X)$, and the identity hom $1_{G(A)}\in\operatorname*{Hom}_{\mathbf{X}%
}(G(A),G(A))$ is associated by the right isomorphism with the
\textit{universal sending het} $e_{A}:G(A)\dashrightarrow A$. Then for any het
$d\in\operatorname*{Het}(X,A)$, the uniquely associated hom $f(d)\in
\operatorname*{Hom}_{\mathbf{A}}(F(X),A)$ gives: $d=f\left(  d\right)  h_{X}$.
And for the same $d$, the uniquely associated $g\left(  d\right)
\in\operatorname*{Hom}_{\mathbf{X}}(X,G(A))$ gives: $d=e_{A}g(d)$. Combining
the two gives the adjunctive square diagram used previously in the text.%

\begin{center}
\includegraphics[
height=2.4882in,
width=5.0394in
]%
{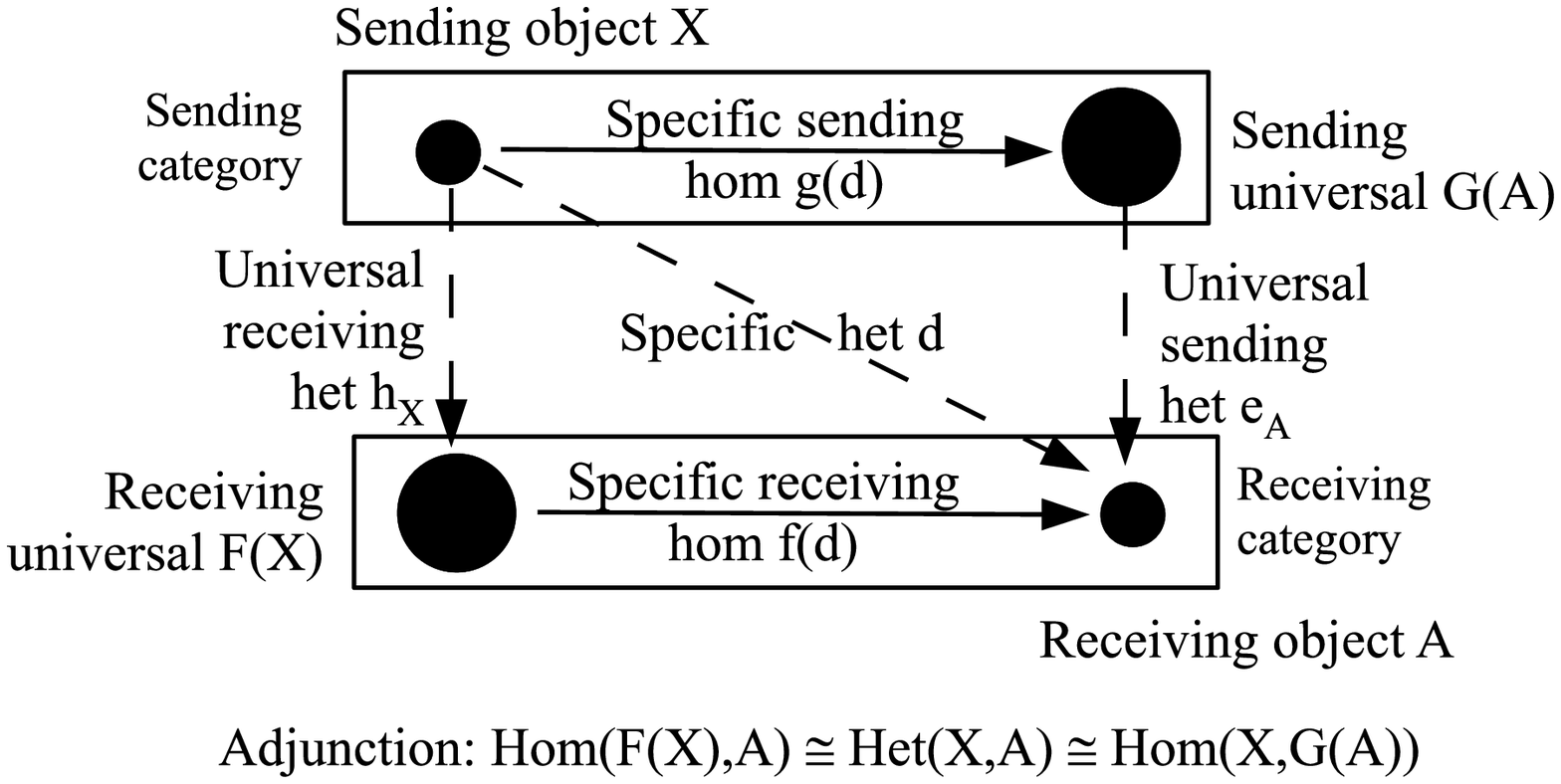}%
\\
Figure 23: Adjunctive square diagram.
\end{center}

Finally, a \textit{brain functor} is a functor $F:\mathbf{X}\rightarrow
\mathbf{A}$ that is a left semiadjunction for $\operatorname*{Het}(X,A)$ and a
right semiadjunction for $\operatorname*{Het}(A,X)$, i.e.,

\begin{center}
$\operatorname*{Hom}_{\mathbf{A}}(F(X),A)\cong\operatorname*{Het}(X,A)$ and
$\operatorname*{Het}(A,X)\cong\operatorname*{Hom}_{A}(A,F(X))$.
\end{center}

\noindent For each $d\in\operatorname*{Het}\left(  X,A\right)  $, there is a
unique hom $f\left(  d\right)  \in\operatorname*{Hom}_{\mathbf{A}}(F(X),A)$ so
that the upper triangular `wing' in the butterfly diagram commutes. For each
$d^{\prime}\in\operatorname*{Het}\left(  A,X\right)  $, there is a unique hom
$g(d^{\prime})\in\operatorname*{Hom}_{A}(A,F(X))$ so that the lower triangular
`wing' commutes.%

\begin{center}
\includegraphics[
height=2.8028in,
width=4.7223in
]%
{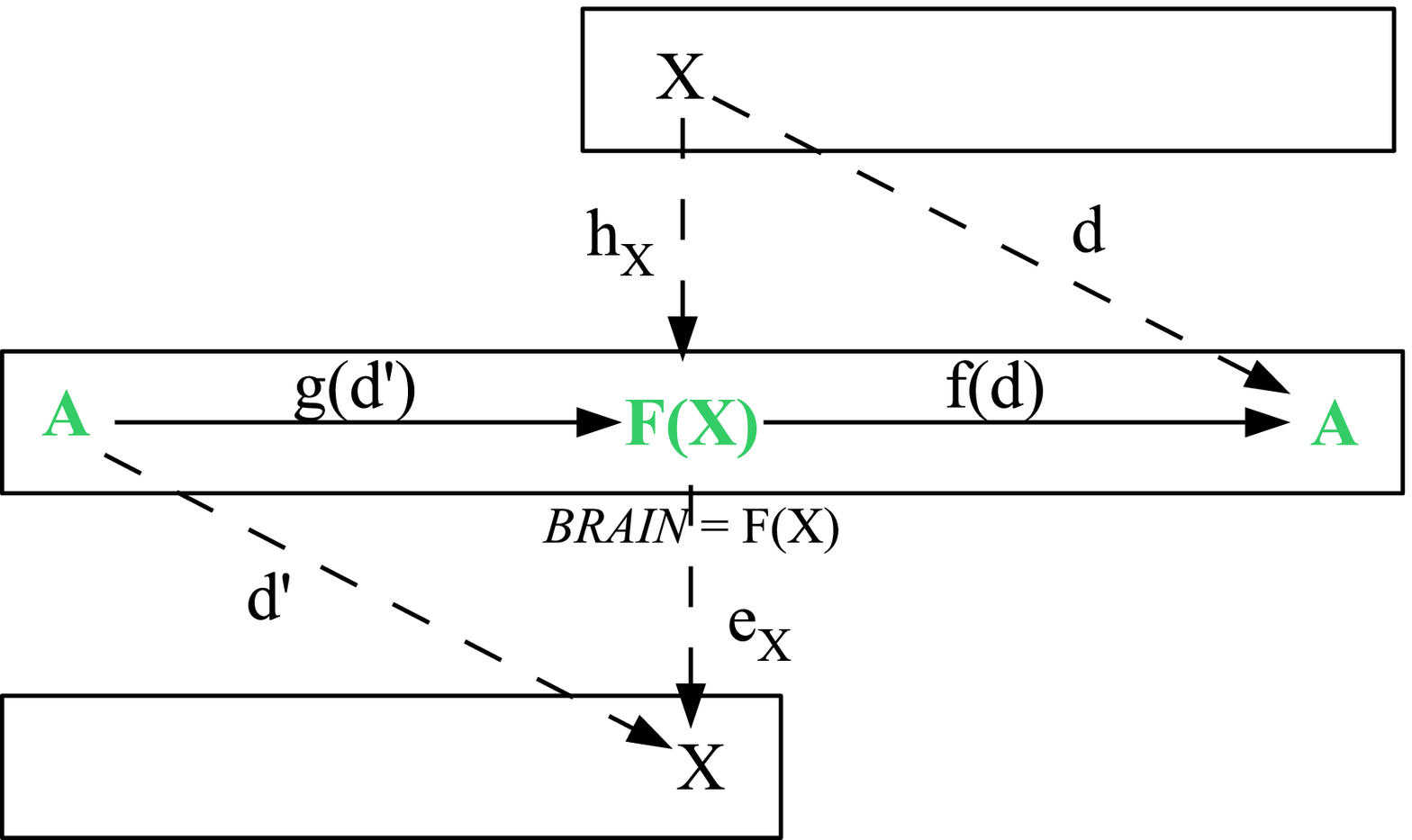}%
\\
Figure 24: Mathematical butterfly diagram for a brain functor.
\end{center}

If a functor $F:\mathbf{X}\rightarrow\mathbf{A}$ has a right adjoint
$G:\mathbf{A}\rightarrow\mathbf{X}$, then:

\begin{center}
$\operatorname*{Hom}_{\mathbf{A}}(F(X),A)\cong\operatorname*{Het}%
(X,A)\cong\operatorname*{Hom}_{\mathbf{X}}(X,G(A))$.
\end{center}

\noindent If the functor $F$ also has a left adjoint $H:\mathbf{A}%
\rightarrow\mathbf{X}$, then:

\begin{center}
$\operatorname*{Hom}_{\mathbf{X}}(H(A),X)\cong\operatorname*{Het}\left(
A,X\right)  \cong\operatorname*{Hom}_{\mathbf{A}}(A,F(X))$.
\end{center}

\noindent Then taking the isomorphisms that do not involve $G$ or $H$ gives
$\operatorname*{Hom}_{\mathbf{A}}(F(X),A)\cong\operatorname*{Het}(X,A)$ and
$\operatorname*{Het}\left(  A,X\right)  \cong\operatorname*{Hom}_{\mathbf{A}%
}(A,F(X))$, i.e., $F$ is a brain functor. Hence all functors $F$ that have
both right and left adjoints are brain functors.

\end{document}